\documentstyle[twoside, 12pt]{article}
\input amssym.def

\newenvironment{demo}[1]%
{\vskip-\lastskip\medskip
  \noindent
  {\em #1.}\enspace
  }%
{\qed\par\medskip
  }

\newcommand{\qed}{
  \strut\hfill
  \mbox{$\Box$}
  }

{\vskip-\lastskip
  \begin{trivlist}
    \item[]
      {\bf Axiom #1 }
      }%
{\end{trivlist}
  }

\newtheorem{theorem}{Theorem}[section]

\newtheorem{corollary}{Corollary}[section]

\newtheorem{lemma}{Lemma}[section]

\newtheorem{remark}{Remark}[section]

\newcommand{\binf}{ b_{\infty} }

\newcommand{\binftwo}{ \tilde{b}_{\infty} }

\newcommand{\Binf}{ {\overline{b}}_{\infty} }

\newcommand{\Binftwo}{ {\overline{b^{\tilde{}}}}_{\infty} }

\newcommand{\cinf}{ c_{\infty} }

\newcommand{\cinfpm}{  {\cinf}_{\pm} }

\newcommand{\Cinf}{ {\overline{c}}_{\infty} }

\newcommand{\dinf}{ d_{\infty} }

\newcommand{\dinfpm}{  {\dinf}_{\pm} }

\newcommand{\Dinf}{ {\overline{d}}_{\infty} }

\newcommand{\C}{ {\Bbb C} }

\newcommand{\Fhalf}{ {\cal F}^{\bigotimes \hf} }

\newcommand{\Fhalfminusl}{ {\cal F}^{\bigotimes -l + \hf} }

\newcommand{\Fl}{ {\cal F}^{\bigotimes l} }

\newcommand{\Fm}{ {\cal F}^{\bigotimes m} }

\newcommand{\Fn}{ {\cal F}^{\bigotimes n} }

\newcommand{\Fmn}{ {\cal F}^{\bigotimes (m + n)} }

\newcommand{\Flhalf}{ {\cal F}^{\bigotimes l + \hf} }

\newcommand{\Flminushalf}{ {\cal F}^{\bigotimes l - \hf} }

\newcommand{\Flpmhalf}{ {\cal F}^{\bigotimes l \pm \hf} }

\newcommand{\Fminushalf}{ {\cal F}^{\bigotimes  - \hf} }

\newcommand{\Fminusl}{ {\cal F}^{\bigotimes  - l} }

\newcommand{\Fminuslhalf}{ {\cal F}^{\bigotimes  - l - \hf} }

\newcommand{\Fminuslpmhalf}{ {\cal F}^{\bigotimes  - l \pm\hf} }

\newcommand{\Fpair}{ ( O(2l), \dinf ) }

\newcommand{\Fpairmtimesn}{ ( O(2m) \times O(2n), \dinf \times \dinf) }

\newcommand{\Fpairbb}{ ( Spin(2l+1), \binf ) }

\newcommand{\Fpairbbtwo}{ ( Spin(2l+1), \binftwo  ) }

\newcommand{\Fpairbd}{ ( O(2l+1), \dinf ) }

\newcommand{\Fpaircc}{ (Sp (2l), \cinf ) }

\newcommand{\Fpairdb}{ ( Pin(2l), \binf ) }

\newcommand{\Fpairdbtwo}{ ( Pin(2l), \binftwo  ) }

\newcommand{\Fpairgl}{ ( GL(l), \hgl ) }

\newcommand{\Fpairospc}{ (Osp (1, 2l), \cinf ) }

\newcommand{\gl}{ {\frak {gl}} }

\newcommand{\glpm}{ \widehat{{ \frak {gl} }}_{\pm} }

\newcommand{\glzero}{ \widehat{{ \frak {gl} }}_0 }

\newcommand{\hf}{ \frac12}

\newcommand{\hgl}{ \widehat{ \frak{gl} } }

\newcommand{\hn}{ \hf + {\Bbb Z}_{+} }

\newcommand{\hz}{ \hf + \Bbb Z }

\newcommand{\hL}{ \widehat{\Lambda} }

\newcommand{\Mpairbc}{ ( O(2l +1), \cinf ) }

\newcommand{\Mpaircd}{ ( Sp(2l), \dinf ) }

\newcommand{\Mpairdc}{ ( O(2l ), \cinf ) }

\newcommand{\Mpairgl}{ ( GL(l), \hgl ) }

\newcommand{\Mpairospd}{ (Osp (1,2l), \dinf) }

\newcommand{\NN}{ z^{-n-1} }

\newcommand{\Od}{ {}^{\frak d} {\cal O} }

\newcommand{\Odf}{ {}^{\frak d} {\cal O}_f }

\newcommand{\Odfm}{ {}^{\frak d} {\cal O}_f^m }

\newcommand{\Odfn}{ {}^{\frak d} {\cal O}_f^n }

\newcommand{\Tr}{ \mbox{Tr } }

\newcommand{\SUM}{ \sum_{n\in \hf + \Bbb Z} }

\newcommand{\vac}{ |0 \rangle }

\newcommand{\W}{ {\cal W}_{1+\infty} }

\newcommand{\Z}{ {\Bbb Z} }

\begin{document}
\title{ 
 Duality in infinite dimensional Fock representations
  }
\author{
  %
  Weiqiang Wang\\
\\{\small Max-Planck Institut f\"ur Mathematik, 53225 Bonn, Germany}\thanks{
On leave from Department of Mathematics, Yale University,
New Haven, CT 06520, USA}
\\{\small E-mail: wqwang@mpim-bonn.mpg.de}
}
\date{}
\maketitle
\begin{abstract}
We construct and study in detail various dual pairs 
acting on some Fock spaces between a finite dimensional 
Lie group and a completed infinite rank affine algebra 
associated to an infinite affine Cartan matrix. 
We give explicit
decompositions of a Fock representation into a direct sum
of irreducible isotypic subspaces with respect to the action of a
dual pair, present explicit formulas for the common highest
weight vectors and calculate the corresponding highest weights. 
We further outline applications of these dual pairs 
to the study of tensor products of modules of such an 
infinite dimensional Lie algebra.
\end{abstract}
\setcounter{section}{-1}
\section{Introduction}
 As simple and elegant as any fundamental concept should be,
the theory of dual pairs of Howe has been very successful in the study of
representation theory of reductive groups
(cf. \cite{H1, H2, KV} and references therein). The idea roughly
goes as follows: assume that there is a maximally commuting
pair of Lie groups/algebras acting natually on a vector
space which in turn by itself usually is a ``minimal'' module of
some other larger groups/algebras. Then one gets a decomposition
of the vector space into a direct sum of 
isotypic subspaces which are irreducible under the joint action
of these two commuting Lie groups/algebras. 
Often time one can obtain all the unitary highest
weight representations of a fixed reductive Lie group by 
varying dual pair partners and minimal modules accordingly.
Among many applications, we mention the applications to branching
laws and to a decomposition of a tensor product of
two modules, cf. \cite{H2}.

It is not too surprising that
a fundamental principle such as dual pairs applies to representation
theory of infinite dimensional Lie algebras and superalgebras as well.
The purpose of this paper is to present a systematic study
of dual pairs between a finite dimensional Lie group and
a completed infinite rank affine algebra of Kac-Moody type associated
to an infinite affine Cartan matrix \cite{DJKM, K} acting on some 
infinite dimensional Fock spaces (namely representations of
infinite dimensional Clifford/Heisenberg algebras in plain language). 
One of our applications we have in mind
will be on the study of quasifinite highest weight
modules of Lie subalgebras of $\W$ in 
\cite{KWY} which generalize the earlier works on
the study of representation theory 
of $\W$ algebra \cite{KR1, FKRW, KR2}.
A duality result involving a particular
dual pair of this sort between a finite dimensional
general linear group and an infinite dimensional general
linear Lie algebra $\hgl$ in a fermionic Fock space 
was first obtained by I. Frenkel \cite{F1} by other methods 
than invoking the principle of dual pairs. 

Dual pairs in infinite dimensional 
Fock spaces are intimately related to the classical dual pairs in 
finite dimensional cases.
However they are not simply direct limits of some classical dual pairs. 
Rather there is some sort of 
semi-infinite twist in the Fock space which, in our opinion, makes
things more interesting. As a consequence,
the infinite dimensional Lie algebras appearing in these
dual pairs have non-trivial central extensions and the irreducible
representations 
appearing in the Fock space decomposition are of highest weight. 
The central charge of these representations corresponds
to the rank of the Lie group in the corresponding dual pair.
To some extend, the existence of dual pairs in our infinite dimensional 
setting ensures the stability properties of classical dual pairs.

It turns out that the Fock spaces on which we have 
constructed various dual pairs are those which
A.~Feingold and I.~Frenkel used to realize the level $\pm 1$
representations of classical affine algebras \cite{FF}. A key observation
is that the Lie algebras of the
Lie groups appearing in our dual pairs are the horizontal 
subalgebras of these classical affine algebras accordingly,  
both untwisted and twisted cases included.

There will be various dual pairs, depending on whether
a Fock space is built out of fermionic or bosonic
fields, whether these fields are indexed by $\Z $ or $\hz $, 
whether there is an additional neutral field or not,
whether the related affine Kac-Moody algebra is untwisted or twisted. 
We list below in two tables the dual pairs which we will 
construct in this paper in infinite dimensional
Fock spaces. The first (resp. second) table lists those dual pairs appearing
in a Fock space of $l$ pairs of fermionic (resp. bosonic ghost)
fields possibly together with 
an extra fermionic or bosonic field.

\begin{table}
 \begin{center}
 \vspace{0.1in}
  \begin{tabular}
 {|l |l |l |l |l |} \hline
           &          &$\Fl$         &$\Flhalf$     &$\Flminushalf$ \\ \hline
 Untwisted& $\hf+\Z$  & $\Fpair  $   & $\Fpairbd $  &               \\ \hline
 Untwisted &  $ \Z$   & $\Fpairdbtwo$& $\Fpairbbtwo$ &              \\ \hline
 Twisted   & $\hf+\Z$ & $\Fpaircc$   &            & $\Fpairospc$    \\ \hline
 Twisted   & $  \Z  $ & $\Fpairdb$   & $\Fpairbb$    &            \\
  \hline
 \end{tabular} 
 \end{center}
 \caption{\protect dual pairs in (mostly) fermionic Fock spaces}
\end{table} 

\begin{table}
 \begin{center}
 \begin{tabular}
 {|l |l |l |l |l |} \hline
           &          &$\Fminusl$    &$\Fhalfminusl$ &$\Fminuslhalf$\\ \hline
 Untwisted& $\hf+\Z$  & $\Mpaircd$   & $\Mpairospd$  &              \\ \hline
           %
                  %
 Twisted   & $\hf+\Z$ & $\Mpairdc$   &               & $\Mpairbc$   \\ \hline
                  %
 \end{tabular} 
 \end{center}
 \caption{\protect dual pairs in (mostly) bosonic Fock spaces}
\end{table} 
 \vspace{0.1in}

Some remarks on notations are in order. We use ${\cal F}^{\bigotimes *}$
to denote various Fock spaces throughout this paper.
The plus or minus sign in the tensor power of a Fock
space means the Fock space is
fermionic or bosonic. The number $l$ in the tensor
power of a Fock space indicates it is a Fock space of
$l$ pairs of fields which are fermions or bosonic ghosts depending
on the sign. $ \pm \hf$ in the tensor
power indicates the appearance of an extra neutral field,
fermionic or bosonic depending on the sign again. For instance,
$\Fl$ denotes the Fock space of $l$ pairs of fermionic fields. 
$\Fhalfminusl$ denotes the Fock space of $l$ pairs of bosonic ghost fields
and a neutral fermionic field. 

The entry $\Z$ or $\hf + \Z$ in the two tables
means the Fourier components are indexed by $\Z$ or $ \hf +\Z$.
Untwisted (resp. twisted) in the two tables means that 
the Lie algebra of the finite dimensional Lie group appearing 
in the corresponding dual pair is the horizontal Lie subalgebra 
of an untwisted (resp. twisted) classical affine algebra
acting on the same Fock space.
For instance, the dual pair $\Fpairbd $ in the first table appears in a row 
starting with ``Untwisted'' and $\hf+\Z$ and in a column starting
with $\Flhalf$. This tells us that there is an untwisted
affine algebra (which is indeed $\widehat{ \frak{so}}(2l +1)$ in this case)
acting on the Fock space $\Flhalf$ of $l$ pairs of fermionic fields 
and a neutral fermionic field which are indexed by $\hf + \Z$. 
The action of $O(2l +1)$ in the dual pair
is given by integrating the action of the horizontal subalgebra
$ \frak{so}(2l +1)$ of $\widehat{ \frak{so}}(2l +1)$.
$\binf$, $\binftwo$, $\cinf$ and $\dinf$ appearing in the
tables are Lie subalgebras of $\hgl$ of $B, C, D$ types
as defined in Section \ref{sec_algebras}, where $\hgl$ is
the central extension of the Lie algebra of infinite size matrices
$(a_{ij})_{i,j \in \Z}$ with finitely many non-zero diagonals.

The dual pair between a general linear group $GL(l)$ and
$\hgl$ is not listed in these two tables. Indeed we have a dual
pair $\Fpairgl$ acting on the Fock spaces $\Fl$ and $\Fminusl$ 
regardless the other constraints
on the tables, such as twisted or untwisted, $\Z$ or $\hf +\Z$.
Dual pair $\Mpairgl$ acting on the Fock space $\Fminusl$ was
treated earlier by Kac and Radul \cite{KR2}, 
where the decomposition of $\Fminusl$ 
into isotypic subspaces was given. However
the calculation of the highest weights in isotypic spaces was performed
there in some indirect way without appealing to
explicit formulas for the highest weight vectors in isotypic spaces. 

In each case listed in the two tables, 
we completely determine the Fock space decomposition
into a direct sum of isotypic subspaces with respect to 
a corresponding dual pair, give explicit formulas
for the highest weight vectors in these isotypic subspaces
and calculate their highest weights. The highest weights can be read off 
from the explicit formulas of these highest weight vectors.
We will see that 
these Fock spaces usually provides natural models on which 
every finite dimensional irreducible representation 
of a finite dimensional Lie group
in the corresponding dual pair appears and every unitary highest weight
module of an infinite dimensional Lie algebra 
in the corresponding dual pair also appears.		

While the theory of classical dual pairs has many 
far-reaching applications \cite{H1, H2}, we confine ourselves
to the discussion of one particular application in the paper.
We establish a number of reciprocity laws associated to
see-saw pairs arising in these Fock spaces.
We construct semisimple tensor categories 
of certain highest weight modules of $\dinf$ 
(resp. $\binf$, $\binftwo$, $\cinf$)
and establish the equivalence of tensor categories with
some suitable tensor categories of representations of
finite dimensional Lie groups with various ranks.
Similar results on tensor categories were
studied earlier by the author
in the case of dual pair $(GL(l), \hgl)$ \cite{W}.

One may have various formulations of our duality results
just as in the case of classical dual pairs.
We do not intend to do so since it is quite clear
how to do them by imitating \cite{H2} and doing that will increase the
length of the paper considerably.

The paper is organized as follows. In Section \ref{sec_algebras}
we review Lie algebra $\hgl$ and its subalgebra 
$\binf, \binftwo, \cinf, \dinf$ of
$B, C, D$ types. In Section \ref{sec_classical} we 
give a parametrization
of finite dimensional irreducible representations of classical
Lie groups and $Spin(2l), Pin(2l),  Osp (1, 2l)$ 
which appear in our Fock space decompositions.
In Section \ref{sec_Fpair} we study the
dual pair actions in the fermionic Fock space $\Fl$.
In Section \ref{sec_Fpairbd} we study the
dual pair action in the Fock space $\Flpmhalf$.
In Section \ref{sec_Bpair} we work out in detail duality
in the bosonic Fock space $\Fminusl$.
In Section \ref{sec_super} we treat dual pairs
in the Fock space $\Fminuslpmhalf$. In Section \ref{sec_reciprocity}
we outline a number of reciprocity laws and a study of
various tensor categories.

Conventions: $\Z $ is the set of integers; $\Z_{+} $ is the
set of non-negative integers; $\Bbb N $ is the set of positive integers; 
$\underline{\Bbb Z}$ is either 
$\Z$ or $\hz$. $\underline{\Bbb Z}_{+} $ is $\hf +{\Z_{+} }$ if
$\underline{\Bbb Z} = \hf + \Z$, and $\Bbb N$ if $\underline{\Bbb Z} = \Z$.
All classical Lie groups and algebras are over the complex field $\C$
unless otherwise specified. An irreducible representation of
a finite dimensional Lie group/algebra always
means to be finite dimensional. Sometimes we will simply
use $\lambda$ to denote the highest weight representation
with highest weight $\lambda$ when no ambiguity may arise.

\noindent{\bf Acknowledgement.}
   I am grateful to R.~Howe for very helpful correspondences and remarks
to C.~H.~Yan for helpful discussions. 
I also thank Max-Planck-Institut for finaicial support and excellent working
environment.
\section{Lie algebras $\hgl$, $\binf$, $\binftwo$, $\cinf$ and $\dinf$ }
\label{sec_algebras}
In this section we review and fix notations on Lie algebras 
$\hgl$ and its various Lie subalgebras of $B, C, D$ types, cf. e.g. \cite{K}.
\subsection{Lie algebra $\hgl$}
  Let us denote by $\gl_f$ the Lie algebra of all matrices
of infinite size
$(a_{ij})_{i,j \in \Z}$ with finitely many non-zero entries.
Denote by $\gl$ the Lie algebra of all matrices
$(a_{ij})_{i,j \in \Z}$ with only finitely many nonzero
diagonals. Obviously $\gl_f$ is a Lie subalgebra of $\gl$.
Putting weight $ E_{ij} = j - i $ defines
a $\Z$--principal gradation 
$\gl = \bigoplus_{j \in \Z} \gl_j$.
Denote by $\hgl = \gl \bigoplus \C C$ the central extension
given by the following $2$--cocycle with values in $\C$ \cite{DJKM}:
\begin{eqnarray}
 C(A, B) = \Tr \,\left(
                     [J, A]B \right)
  \label{eq_cocy}
\end{eqnarray}
where $J = \sum_{ j \leq 0} E_{ii}$. The $\Z$--gradation
of Lie algebra $\gl$ extends to $\hgl$ by putting weight $C = 0$. 
In particular, we have a triangular decomposition
$$\hgl = \widehat{\gl}_{+}  
          \bigoplus \widehat{\gl}_{0} \bigoplus \widehat{\gl}_{-}  $$
where 
$$ \glpm = \bigoplus_{ j \in \Bbb N}
\widehat{\gl}_{\pm j}, \quad 
\widehat{\gl}_{0} = \gl_0 \oplus \C C.$$
Denote by $E_{ij}$ the infinite matrix with $1$ at $(i, j)$
place and $0$ elsewhere. Denote by $\epsilon_i$ the linear
function on ${\frak {gl} }_0$, s.t. 
$\epsilon_i (E_{jj} ) = \delta_{ij} (i, j \in \Z ).$
Then the root system of $\gl$ is 
$\Delta = \{ \epsilon_i - \epsilon_j, (i, j \in \Z, i \neq j) \}.$
The compact anti-involution $\omega$ is defined as
$\omega (E_{ij}) = E_{ji}$.

Given $c \in \Bbb C$ and $\Lambda \in {\gl}_{0}^{*}$, we let
\begin{eqnarray*}
 \lambda^a_i & = & \Lambda (E_{ii}), \quad i \in \Z,    \\
 H^a_i     & = & E_{ii} - E_{i+1, i+1} + \delta_{i,0} C,  \\
  h^a_i & = & \Lambda( H^a_i ) = 
   \lambda_i - \lambda_{i+1} + \delta_{i,0} c.
 \label{eq_11}
\end{eqnarray*}
The superscript $a$ here denotes $\hgl$ which is of A type.

Denote by $L(\hgl; \Lambda, c)$ (or simply $L(\hgl; \Lambda)$ 
when the central charge is obvious from the text)
the highest weight $\hgl$--module with
highest weight $\Lambda$ and central
charge $c$. Easy to see that $L(\hgl; \Lambda, c)$ is {\em quasifinite}
(namely having finite dimensional graded subspaces
according to the principal gradation of $\hgl$) if
and only if all but finitely many $h_i, i \in \Z$ are zero.
A quasifinite representation of $\hgl$ is {\em unitary}
if an Hermitian form naturally defined with respect to $\omega$
is positive definite. 

Define $\Lambda_j^a \in \gl_0^* , j \in \Z$ as follows:
\begin{equation}
\Lambda_j^a ( E_{ii} ) =    
   \left\{
      \everymath{\displaystyle}
      \begin{array}{lll}
        1, & \mbox{for}\quad 0 < i \leq j \\
        -1, & \mbox{for}\quad j < i \leq 0  \\
        0, & \mbox{otherwise.}
      \end{array}
    \right. \\
  \label{eq_122}
\end{equation}
Define $\hL_0^a \in \hgl_0^{*}$ by 
$$ \hL_0^a (C) = 1, \quad
  \hL_0^a (E_{ii}) = 0 \mbox{ for all } i \in \Z $$
and extend $\Lambda_j^a\/$ from $ \gl_0^{*}\/$ to $\hgl_0^{*}$ by letting
$\Lambda_j^a (C) = 0.$ Then 
$$\hL_j^a = \Lambda_j^a + \hL_0^a, \quad j \in \Z  $$ 
are the fundamental weights,
i.e. $\hL_j^a ( H_i ) = \delta_{ij}.$

It is not difficult to prove that 
$L(\hgl; \Lambda, c)$ is unitary if and only if
$\Lambda =  \hL^a_{m_1} + \ldots + \hL^a_{m_k}, c = k \in \Z_{+}$
by using a method due to Garland \cite{G}.
\subsection{Lie algebra $\dinf$}
  \label{subsec_dinf}
   Now consider the vector space $\C [t, t^{-1}]$ and take a basis
$v_i = t^i, i \in \Z$. The Lie algebra $\gl$ acts on this vector
space naturally, namely $E_{ij} v_k = \delta_{jk}v_i$. We denote
by $\Dinf$ the Lie subalgebra of $\gl$ preserving the following
symmetric bilinear form (cf. \cite{K}):
$$ D(v_i, v_j) = \delta_{i, 1-j}, \quad i, j \in \Z.$$
Namely we have 
$$\Dinf = \{ g \in \gl \mid D( a(u), v) + D(u, a(v)) = 0 \}
 = \{ (a_{ij})_{i,j \in \Z} \in \gl \mid a_{ij} = -a_{1-j,1-i} \}.$$
Denote by $\dinf = \Dinf \bigoplus \C C $ the central extension
given by the $2$-cocycle (\ref{eq_cocy}) restricted to $\Dinf$.
Then $\dinf$ has a natural triangular decomposition induced
from $\hgl$:
$$\dinf = {\dinf}_{+} \bigoplus {\dinf}_0 \bigoplus {\dinf}_{-} $$
where 
$\dinfpm = \dinf \cap \glpm $ and $ {\dinf}_0 = {\dinf} \cap \glzero .$

The set of simple coroots of $\dinf$, denoted by $ \Pi^\vee $
can be described as follows:
\begin{eqnarray*}
      \begin{array}{rcl}
      \lefteqn{
        \Pi^\vee = \{ \alpha_0^\vee = E_{00} + E_{-1,-1}
              -E_{2, 2} -E_{1, 1} + C,
              }                   \\
      & & \alpha_i^\vee = E_{i, i} + E_{-i,-i}
        - E_{i+1,i+1} - E_{1-i,1-i}, i \in \Bbb N \}.
      \end{array}
\end{eqnarray*}
               %
               %
               %
               %
               %
               %

Given $\Lambda \in {\dinf}_0^* $, we let
\begin{eqnarray*}
\lambda^d_i & = & \Lambda (E_{ii} - E_{1-i, 1-i}) \quad (i \in \Bbb N), \\
 H^d_i & = & E_{ii} + E_{-i, -i} - E_{i+1, i+1} - E_{-i+1, -i+1}
                                \quad (i \in \Bbb N),   \\
  h^d_i & = & \Lambda( H^d_i ) = 
   \lambda_i  - \lambda_{i+1}  \quad (i \in \Bbb N), \\
 H^d_0 & = & E_{0,0} + E_{-1,-1} -E_{2,2} -E_{1,1} + 2C,   \\
 c & = & \hf (h^d_0 + h^d_1) + \sum_{i \geq 2} h^d_i.
 \label{eq_dcartan} 
\end{eqnarray*}
The superscript $d$ denotes $\dinf$ which is of D type.
Then we denote by $\hL^d_i $ the $i$-th fundamental weight
of $\dinf$, namely $\hL^d_i (H^d_j ) = \delta_{ij}.$
Denote by $L(\dinf; \Lambda, c)$
(or $L(\dinf; \Lambda)$ when the central charge is obvious)
the highest weight $\dinf$--module with
highest weight $\Lambda$ and central charge $c$. 
Such a representation is {\em unitary}
if an Hermitian form naturally defined with respect to
the compact anti-involution $\omega$ on $\hgl$ when restricted to
$\dinf$ is positive definite. It is not difficult to show that 
$L(\dinf, \Lambda)$ is unitary if and only if
$\Lambda =  \hL^d_{m_1} + \ldots + \hL^d_{m_k}.$

\subsection{Lie algebras $\binf$ and $\binftwo$}
   \label{subsec_binf}
 Let us consider the following symmetric bilinear forms
$B(v_i, v_j ) = ( -1)^i \delta_{i, -j}$ and 
$\widetilde{B}(v_i, v_j ) = \delta_{i, -j}, i, j \in \Z.$ Denote by
$\Binf$ (resp. $\Binftwo$) the Lie subalgebra of $\gl$ which preserves the
bilinear form $B$ (resp. $\widetilde{B}$), namely we have 
\begin{eqnarray*}
 \Binf & = & \{ g \in \gl \mid B( a(u), v) + B (u, a(v)) = 0 \}  \\
        & = &\{ (a_{ij})_{i,j \in \Z} \in \gl
               \mid a_{ij} = - ( -1)^{i+j}a_{-j,-i} \}.  \\
   \Binftwo & = & \{ g \in \gl \mid
                \widetilde{B}( a(u), v) + \widetilde{B}(u, a(v)) = 0 \}  \\
        & = & \{ (a_{ij})_{i,j \in \Z} \in \gl \mid a_{ij} = - a_{-j,-i} \}.
\end{eqnarray*}
Denote by $\binf = \Binf \bigoplus \C C$ 
(resp. $\binftwo = \Binftwo \bigoplus \C C$)
the central extension of $\Binf$ (resp. $\Binftwo$) 
given by the $2$-cocycle
(\ref{eq_cocy}) restricted to $\Binf$ (resp. $\Binftwo$). 
Clearly these two Lie algebras $\binf$ and $\binftwo$
share the same Cartan subalgebras. Then $\binf$
(resp. $\binftwo$) inherits from $\hgl$ a natural triangular decomposition:
\begin{eqnarray*}
 \binf = {\binf}_{+} \bigoplus {\binf}_0 \bigoplus {\binf}_{-},   \quad
 \binftwo = \tilde{b}_{\infty,+}
               \bigoplus {\binf}_0 \bigoplus \tilde{b}_{\infty,-}
\end{eqnarray*}
where ${b}_{\infty,\pm} = \binf \cap \glpm $,
$ {\binf}_0 = {\binf} \cap \glzero ,$
and  $\tilde{b}_{\infty,\pm} = \binftwo \cap \glpm $.

The set of simple coroots of $\binf$ (which is the same as
that of $\binftwo$), denoted by $ \Pi^\vee $
can be described as follows:
\begin{eqnarray*}
      \begin{array}{rcl}
      \lefteqn{
          \Pi^\vee = \{ \alpha_0^\vee =2( E_{-1,-1} - E_{1,1}) + C,  
              }             \\
       & & \alpha_i^\vee = E_{i, i} + E_{-i-1,-i-1}
              - E_{i+1,i+1} - E_{-i,-i}, i \in \Bbb N \}.
      \end{array}
\end{eqnarray*}
               %
                %
                %
               %
               %
               %

Given $\Lambda \in {\binf}_0^* $, we let
\begin{eqnarray*}
 \lambda^b_i & = & \Lambda (E_{ii} - E_{-i, -i}),    \\
  H^b_i & = & E_{ii} + E_{-i-1, -i-1} - E_{i+1, i+1} - E_{-i, -i},    \\
   H^b_0 & = & 2 (E_{-1,-1} -E_{1,1}) + 2C,            \\
  h^b_i & = & \Lambda( H^b_i ) = \lambda_i - \lambda_{i+1},
  \quad i \in \Bbb N,                              \\
  c & = & \hf h^b_0 +  \sum_{i \geq 1} h^b_i.
\end{eqnarray*}
The superscript $b$ here denotes $\binf$ and $\binftwo$ which are of B type.
Denote by $\hL^b_i $ the $i$-th fundamental weight
of $\binf$ as well as $\binftwo$, namely
$\hL^b_i (H^b_j ) = \delta_{ij}.$
Denote by $L(\binf; \Lambda)$ (resp. $L(\binftwo; \Lambda)$)
the highest weight module over $\binf$ (resp. $\binftwo$) with
highest weight $\Lambda$ and central charge $c$. 
Such a representation is {\em unitary}
if an Hermitian form naturally defined with respect to
the compact anti-involution $\omega$ on $\hgl$ when restricted to
$\binf$ (resp. $\binftwo$) is positive definite.
It is not difficult to show that 
$L(\binf, \Lambda)$ (resp. $L(\binftwo; \Lambda)$)
is unitary if and only if
$\Lambda = \hL^b_{m_1} + \ldots + \hL^b_{m_k}.$
\subsection{Lie algebra $\cinf$}
Let us consider the skew-symmetric bilinear forms
$C(v_i, v_j ) = ( -1)^i \delta_{i, -j +1}, i, j \in \Z.$ Denote by
$\Cinf$ the Lie subalgebra of $\gl$ which preserves the
bilinear form $C$, namely we have 
\begin{eqnarray*}
 \Cinf & = & \{ g \in \gl \mid
                C( a(u), v) + {C}(u, a(v)) = 0 \}      \\
       & = & \{ (a_{ij})_{i,j \in \Z} \in \gl\mid 
                    a_{ij} = - (-1)^{i+j}a_{1-j,1-i} \}.
\end{eqnarray*}
Denote by $\cinf = \Cinf \bigoplus \C C$ 
the central extension of $\Cinf$ 
given by the $2$-cocycle
(\ref{eq_cocy}) restricted to $\Cinf$. 
Then $\cinf$ inherits from $\hgl$ a natural triangular decomposition:
\begin{eqnarray*}
 \cinf = {\cinf}_{+} \bigoplus {\cinf}_0 \bigoplus {\cinf}_{-} 
\end{eqnarray*}
where $\cinfpm = \cinf \cap \glpm $, $ {\cinf}_0 = {\cinf} \cap \glzero $.

The set of simple coroots of $\cinf$, denoted by $ \Pi^\vee $
can be described as follows:
\begin{eqnarray*}
      \begin{array}{rcl}
      \lefteqn{
          \Pi^\vee = \{ \alpha_0^\vee =E_{0,0} - E_{1,1} + C,  
              }             \\
       & & \alpha_i^\vee = E_{i, i} + E_{-i,-i}
              - E_{i+1,i+1} - E_{1-i,1-i}, i=1,2,... \}.
      \end{array}
\end{eqnarray*}
Given $\Lambda \in {\cinf}_0^* $, and $c \in \C$, we let
\begin{eqnarray*}
 \lambda_i & = & \Lambda (E_{ii} - E_{1-i, 1-i}),    \\
  H^c_i & = & E_{ii} + E_{-i, -i} - E_{i+1, i+1} - E_{1-i, 1-i},    \\
   H^c_0 & = & E_{0,0} -E_{1,1} + C,            \\
  h^c_i & = & \Lambda( H^c_i ) = \lambda_i - \lambda_{i+1},
  \quad i \in \Bbb N,                              \\
   c & = & \sum_{i \geq 0} h^c_i.
\end{eqnarray*}
The superscript $c$ here denotes $\cinf$ which is of C type.
Then we denote by $\hL^c_i $ the $i$-th fundamental weight
of $\cinf$, namely
$\hL^c_i (H^c_j ) = \delta_{ij}.$
Denote by $L(\cinf; \Lambda, c)$
the highest weight module over $\cinf$ with
highest weight $\Lambda$ and central charge $c$. 
Such a representation is {\em unitary}
if an Hermitian form naturally defined with respect to
the compact anti-involution $\omega$ on $\hgl$ when restricted to
$\cinf$ is positive definite.
It is not difficult to show that 
$L(\cinf, \Lambda)$ is unitary if and only if
$\Lambda =  \hL^c_{m_1} + \ldots + \hL^c_{m_k}.$
\section{Parametrization of irreducible representations
  of classical groups}
    \label{sec_classical}
In this section we give a parametrization of irreducible modules
of finite dimensional Lie groups appearing in our Fock space
decompositions in later sections. See \cite{BtD}
for more detail on $Spin(n), Pin (n)$ and other classical groups.
\subsection{$O(2l)$}
  \label{subsec_evenorth}
   We define 
$O(2l) = \{ g \in GL(2l); {}^t g J g = J \}$ with $J$ equal to
            
 \begin{eqnarray*}
   \left[ \begin{array}{cc}
     0      & I_l       \\
   I_l      & 0   
    \end{array} \right].
 \end{eqnarray*}
Lie group $GL(l) $ can be identified as a subgroup of $O(2l)$
consisting of matrices of the form
 \begin{eqnarray*}
  \left[ \begin{array}{cc}
    g      & 0      \\
    0      & {}^t g^{-1}
         \end{array} \right]
   \end{eqnarray*}
where $g$ is an $l \times l$ non-singular matrix. Here and
below ${}^t g$ denotes the transpose of the matrix $g$.
Lie algebra ${\frak {so}}(2l)$ of $SO(2l)$ consists of matrices of the form
  \begin{eqnarray}
  \left[ \begin{array}{cc}
    \alpha      & \beta       \\
    \gamma      & -{}^t{\alpha} 
         \end{array} \right]
  \label{eq_matr}
  \end{eqnarray}
where $\alpha, \beta, \gamma$ are $l \times l$ matrices and
$ \beta, \gamma$ are skew-symmetric. Lie algebra 
${\frak {gl}} (l)$ is identified with the subalgebra of 
${\frak {so}}(2l)$ consisting
of matrices of the form (\ref{eq_matr}) by putting $\beta$ and 
$\gamma$ to be $0$.
We take a Borel subalgebra ${\frak b} ( {\frak {so}}(2l) )$
of ${\frak {so}}(2l)$ to be the intersection of
${\frak {so}}(2l) $ with the set of 
upper triangular matrices of ${\frak {gl}}(2l)$
and a Cartan subalgebra ${\frak h} ({\frak {so}}(2l) )$ consisting of 
diagonal matrices $diag (t_1, \ldots, t_l, -t_1, \ldots, - t_l ),$
$t_i \in \C.$ The subalgebra ${\frak {gl}}(l)$ 
of ${\frak {so}}(2l) $ share the same
Cartan subalgebra with ${\frak {so}}(2l) $.
         
An irreducible representation of $GL(l)$ is parametrized by its highest
weight with respect to the chosen Cartan subalgebra by
$$\Sigma(A)
\equiv  \{(m_1, m_2, \ldots, m_l ),\quad
m_1 \geq m_2 \geq \ldots \geq  m_l ,
 m_i \in \Z \}. $$
An irreducible representation of $SO(2l)$ is parametrized by its highest
weight with respect to ${\frak h} ({\frak {so}}(2l))$ in 
$ \{ (m_1, m_2, \ldots, m_l ),\quad
 m_1 \geq m_2 \geq \ldots \geq  m_{l-1} \geq | m_l |,
 m_i \in \Z \}$.

$O(2l)$ is a semi-direct product of $SO(2l)$ by $\Z_2$.
Denote by $\tau$ the $2l \times 2l$ matrix 
\begin{eqnarray}
\left[ \begin{array}{cc}
A & B \\
B & A
\end{array} \right]
  \label{eq_tau}
\end{eqnarray}
with 
$ A = diag (1, \ldots, 1, 0), B = diag (0, \ldots, 0, 1). $
Then $\tau \in O(2l) - SO(2l)$ normalizes the Borel subalgebra 
$\frak b$. If $\lambda$ is a representation of $SO(2l)$ of
highest weight $(m_1, m_2, \ldots, m_l )$,
then $\tau . {\lambda}$ has highest weight 
$(m_1, m_2, \ldots, - m_l )$. It follows that
the induced representation of $(m_1, m_2, \ldots, m_l )$ ($m_l \neq 0$)
to $O(2l)$ is irreducible and its restriction to
$SO(2l)$ is a sum of 
$(m_1, m_2, \ldots, m_l )$ and
$(m_1, m_2, \ldots, - m_l )$.
We denote this irreducible representation $\lambda$ of $O(2l)$ by
$(m_1, m_2, \ldots, \overline{m}_l )$, where $m_l > 0$.
If $ m_l = 0$, the representation
$\lambda = (m_1, m_2, \ldots, m_{l-1}, 0 )$ extends to two
different representations of $O(2l)$, denoted by
$\lambda$ and $\lambda \bigotimes det$,
where $det$ is the $1$-dimensional non-trivial representation of
$O(2l)$.
We denote 
\begin{eqnarray*}
  \Sigma(D) =
   & & \left\{ (m_1, m_2, \ldots, \overline{m}_l ) \mid
       m_1 \geq m_2 \geq \ldots \geq m_l > 0, m_i \in \Z;  \right.
                               \\
   & &   (m_1, m_2, \ldots, m_{l-1}, 0 ) \bigotimes det, \\
   & &   \left.   (m_1, m_2, \ldots, m_{l-1}, 0 ) \mid
       m_1 \geq m_2 \geq \ldots \geq  m_{l-1} \geq 0, m_i \in \Z \right\}. 
\end{eqnarray*}
\subsection{$O(2l+1)$}
     We take $O (2l+1)$ to be the form
\begin{eqnarray}
O (2l+1) & = & \{ g \in GL(2l+1); {}^t g J g = J \}, \\
  SO (2l+1) & = & \{ g \in O (2l+1) ; det \; g =1 \},
\end{eqnarray}
where $J$ is the following $(2l+1) \times (2l+1)$ matrix
  \begin{eqnarray*}
  \left[ \begin{array}{ccc}
    0      & I_l   & 0 \\
    I_l    & 0     & 0   \\
    0      & 0     & 1
         \end{array} \right].
  \end{eqnarray*}
Denote by ${\frak {gl}}(2l +1)$ the Lie algebra of $GL(2l +1)$.
The Lie algebra ${\frak {so}}(2l+1)$ is 
the Lie subalgebra of ${\frak {gl}}(2l +1)$ consisting 
of $(2l+1)\times (2l+1)$ matrices
of the form
  \begin{eqnarray}
  \left[ \begin{array}{ccc}
    \alpha      & \beta         & \delta \\
    \gamma      & -{}^t{\alpha} & h      \\
     -{}^t h    & -{}^t{\delta} & 0
         \end{array} \right]
   \label{eq_form}
  \end{eqnarray}
where $\alpha, \beta, \gamma$ are $l \times l$ matrices
and $ \beta, \gamma$ skew-symmetric. 
The Borel subalgebra ${\frak b} ( {\frak {so} } (2l+1) ) $ 
consists of matrices of the form (\ref{eq_form}) by putting 
$\gamma $ and $h$ to be $0$ and $\alpha$ to be upper triangular.
The Cartan subalgebra ${\frak h} ( {\frak {so} } (2l+1) ) $
consists of diagonal matrices of the form 
$ diag ( t_1, \ldots, t_l; -t_1 \ldots -t_l; 0 ),$ $ t_i \in \C $.
An irreducible module of $SO(2l+1)$ is
parametrized by its highest weight
$(m_1, \ldots, m_l ), m_1 \geq \ldots \geq m_l \geq 0, m_i \in \Z.$

  It is well known that
$ O (2l+1)$ is isomorphic to the direct product $ SO (2l+1) \times \Z_2$
by sending the minus identity matrix to $-1 \in \Z_2 = \{ \pm 1 \}.$
Denote by $det$ the non-trivial one-dimensional representation of $O(2l+1)$.
An representation $\lambda$ of $SO(2l+1)$ extends
to two different representations $\lambda$ and $\lambda \bigotimes det$
of $O(2l+1)$, and all irreducible representations of
$O(2l+1)$ is obtained this way.
Then we can parametrize irreducible 
representatiions of $O(2l+1)$ by 
$ (m_1, \ldots, m_l)$ and $ (m_1, \ldots, m_l) \bigotimes det $.
We denote
  $$ \Sigma(B) = 
    \left\{ (m_1, \ldots, m_l), (m_1, \ldots, m_l) \bigotimes det  \mid
       m_1 \geq \ldots \geq m_l \geq 0, m_i \in \Z \right \}.   $$
\subsection{$Spin (n)$ and $Pin (n)$}
  The Pin group $Pin (n)$ is the double covering group
of $ O(n)$, namely we have
$$ 
1 \longrightarrow \Z_2  \longrightarrow Pin (n)
   \longrightarrow O (n) \longrightarrow 1.       
$$
We then define the Spin group $Spin (n)$ to be the inverse
image of $SO (n)$ under the projection from
$ Pin (n)$ to $ O (n) $. Then we have the following exact
sequence of Lie groups:
$$ 
1 \longrightarrow \Z_2  \longrightarrow Spin (n)
   \longrightarrow SO (n) \longrightarrow 1.       
$$
{\bf Case $ {\bf n = 2l}$.}
Denote $ {\bf 1}_l = (1, 1, \ldots, 1) \in \Z^l $ and 
${\bar{\bf 1} }_l  = (1, 1, \ldots, 1, -1) \in \Z^l $.
An irreducible representation of $Spin (2l)$ 
which does not factor to $SO (2l)$ 
 is an irreducible representation of ${\frak {so}} (2l)$
parametrized by its highest weight 
\begin{eqnarray}
  \lambda = \hf {\bf 1}_l + (m_1, m_2, \ldots, m_l)    \label{eq_wtplus}
\end{eqnarray}
or 
\begin{eqnarray}
\lambda = \hf {\bar{\bf 1} }_l + (m_1, m_2, \ldots, -m_l)  \label{eq_wtminus}
\end{eqnarray}
where $m_1 \geq \ldots \geq m_l \geq 0, m_i \in \Z.$

The Pin group $Pin (2l)$ is not connected. There are two possibilities.
First, an irreducible
representation of $Pin (2l)$ factors to that of $O(2l)$,
then we can use the parametrization of irreducible representations
of $O(2l)$ to parametrize these representations of $Pin (2l)$.

Secondly, an irreducible
representation of $Pin (2l)$ is induced from an 
irreducible representation of $Spin (2l)$ with highest
weight of (\ref{eq_wtplus}) or (\ref{eq_wtminus}).
When restricted to $Spin (2l)$, 
it will decompose into a sum of the two irreducible representations
of highest weights (\ref{eq_wtplus}) and (\ref{eq_wtminus}). 
We will use 
$\lambda = \hf |{\bf 1}_l | + (m_1, m_2, \ldots, \overline{m}_l), m_l \geq 0$
to denote this irreducible representation of $Pin (2l)$. 
Denote by 
$$ \Sigma(Pin) = \{ \hf |{\bf 1}_l | + (m_1, m_2, \ldots, \overline{m}_l), 
   m_1 \geq \ldots \geq m_l \geq 0, m_i \in \Z \}.  $$

{\bf Case ${\bf n = 2l +1}$.}
   An irreducible representation of $Spin (2l +1)$ 
which does not factor to $SO (2l +1)$ 
is an irreducible representation of ${\frak {so}} (2l +1)$
parametrized by its highest weight 
\begin{eqnarray}
\lambda = \hf {\bf 1}_l + (m_1, m_2, \ldots, m_l),  
 \quad m_1 \geq \ldots \geq m_l \geq 0.
  \label{eq_spinodd}
\end{eqnarray}
 
Denote by 
$$ \Sigma(PB) = \left\{ \hf {\bf 1}_l  + (m_1, m_2, \ldots, m_l ) \mid
   m_1 \geq \ldots \geq m_l \geq 0, m_i \in \Z \right \}.   $$
\subsection{$Osp (1, 2l)$ and $Sp (2l)$}
   Denote by $\C^{1|2l}$ the $\Z_2$-graded vector space with
$\C$ as the even subspace and $\C^{2l}$ as the odd subspace.
We denote by $GL(1,2l)$ the general linear Lie supergroup on
the superspace $\C^{1|2l}$, and its Lie superalgebra by
${\frak {gl}} (1, 2l)$. Consider the following supersymmetric
bilinear form
  \begin{eqnarray*}
  \left[ \begin{array}{ccc}
    1     & 0        & 0 \\
    0     & 0        & I_l     \\
    0     & - I_l    & 0
         \end{array} \right]
  \end{eqnarray*}
which is symmetric on the even subspace $\C$ and skew-symmetric
on the odd subspace $\C^{2l}$. 
Define $Osp (1, 2l)$ to be the Lie sub-supergroup of $GL(1, 2l)$
which preserves the above bilinear form. Its Lie superalgebra
${\frak {osp}} (1, 2l) $ 
consists of $(2l +1) \times (2l +1) $ matrices of the following form:
  \begin{eqnarray}
  \left[ \begin{array}{ccc}
    0        & x        & y \\
    {}^t y   & a        & b   \\
    {}^t x   & c    & -{}^t a
         \end{array} \right]
   \label{eq_matrix}
  \end{eqnarray}
where $a, b, c$ are $l \times l$ matrices, $b, c$ are symmetric.

We define ${\frak b} ({\frak {osp}}(1, 2l))$ to be the Borel subalgebra
consisting of matrices of the form (\ref{eq_matrix}) with
$x$ and $c$ to be zero and $a$ to be upper triangular.
The Cartan subalgebra ${\frak h} ({\frak {osp}} (1, 2l) )$
consists of diagonal matrices
$diag \; (0; t_1, \ldots, t_l, -t_1, \ldots, -t_l )$. 

Lie supergroup $Osp (1, 2l)$ shares complete reducibility property
of an ordianry Lie group which makes it distinguished from other
Lie supergroups. So we will treat $Osp (1, 2l)$ just as a Lie group.
An irreducible representation of $Osp (1, 2l)$ can be 
parametrized by its highest weight with respect to 
the Borel and Cartan subalgebra chosen above in
$$ \Sigma (Osp) =
 \{ (m_1, m_2, \ldots, m_l ) \mid 
   \quad m_1 \geq m_2 \geq \ldots \geq m_l \geq 0 \}. $$

The symplectic group $Sp(2l)$ can be taken as the subgroup of
$Osp (1, 2l)$ consisting of matrices in $Osp (1, 2l)$ 
with 1 at the left-upper corner
and 0 at all other entries in the first row and column. We also similarly
define ${\frak b} ( {\frak {sp}} (2l) )$ and
${\frak h} ({\frak {sp}} (2l) )$ to be the Borel and Cartan subalgebra of
${\frak {sp}} (2l)$ consisting of those matrices
with zeros in the first row and column 
in ${\frak b} ({\frak {osp}}(1, 2l))$ and ${\frak h} ({\frak {osp}} (1, 2l) )$
respectively.
An irreducible representation of $Sp (2l)$ can be
parametrized by its highest weight with respect to 
the chosen Borel and Cartan subalgebras by
$$ \Sigma (C) =
 \{ (m_1, m_2, \ldots, m_l ),
     \quad m_1 \geq m_2 \geq \ldots \geq m_l \geq 0  \}. $$
\section{Duality in the fermionic Fock space $\Fl$}
 \label{sec_Fpair} 
\setcounter{subsection}{-1}
\subsection{Fock space $\Fl$}  
     Let us take a pair of fermionic fields
$$ 
 \psi^{+}(z) = \sum_{n \in \underline{\Z} } \psi^{+}_n z^{ -n -\hf + \epsilon},
  \quad \psi^{-}(z) = \sum_{n \in \underline{\Z} }
                           \psi^{-}_n z^{ -n -\hf + \epsilon}, 
  \quad \underline{\Z} = \hf + \Z \mbox{ or } \Z
$$
with 
                 %
                 %
the following anti-commutation relations
$$ [\psi^{+} _m, \psi^{-}_n ]_{+} 
= \psi^{+}_m \psi^{-}_n + \psi^{-}_n \psi^{+}_m = \delta_{m+n, 0}.$$
We take the convention here and below
that $\epsilon = 0$ if $\underline{\Z} = \hz;$ and 
$\epsilon = \hf$ if $ \underline{\Z} = \Z.$
Denote by $\cal F$ the Fock space of the fields 
$\psi^{-}(z)$ and $\psi^{+}(z)$, 
generated by the vacuum $\vac$, satisfying
\begin{eqnarray*}
  \psi^{+}_n \vac = \psi^{-}_n \vac = 0 \quad ( n \in \hn), &&\mbox{ when }
   \underline{\Z} = \hz;                       \\
 \psi^{+}_n \vac = \psi^{-}_{n+1} \vac = 0 \quad (n \in \Z_{+} ),
   && \mbox{ when } \underline{\Z} = \Z. 
\end{eqnarray*}

Now we take $l$ pairs of fermionic fields, 
$\psi^{\pm,p} (z) \;( p = 1, \dots, l)$
and consider the corresponding Fock space $\Fl$.

Introduce the following generating series
\begin{eqnarray}
  E(z,w) & \equiv & \sum_{i,j \in \Z} E_{ij} z^{i-1 + 2\epsilon }
     w^{-j  } 
   = \sum_{p =1}^l :\psi^{+,p} (z) \psi^{-,p} (w):,  \label{eq_genelin} \\
  e^{pq} (z) & \equiv & \sum_{n \in \Z} e^{pq} (n) z^{ -n-1 + 2\epsilon }
    =  : \psi^{-,p}(z) \psi^{-,q} (z): 
    \quad (p \neq q ) ,  \label{eq_orthaffine}    \\ 
  e_{**}^{pq} (z) & \equiv & \sum_{n \in \Z}
                            e_{**}^{pq}(n)z^{ -n-1 + 2\epsilon }
    =  : \psi^{+,p}(z) \psi^{+,q} (z):
    \quad (p \neq q ) , \label{eq_orth}    \\
  e_*^{pq} (z) & \equiv & \sum_{n \in \Z} e_*^{pq} (n) z^{ -n-1 + 2\epsilon }
    =  : \psi^{+,p}(z) \psi^{-,q} (z): +  \delta_{p,q} \epsilon 
         \label{eq_orthaffinestar}  
\end{eqnarray}
where $p,q = 1, \dots, l$, and the normal ordering $::$ means that 
the operators annihilating $\vac$ are moved to the right and 
multiplied by $-1$.

It is well known that the operators $E_{ij} \; (i,j \in \Z)$
generate a representation in $\Fl$ of
the Lie algebra $\hgl$ with central charge $l$;
the operators $e^{pq}(n), e_*^{pq}(n), 
e_{**}^{pq}(n)$, $p,q = 1, \dots, l, n \in \Z$ 
form a representation of the affine algebra
$\widehat{so}(2l)$ with central charge $ 1$  \cite{F1, F, KP}. 
Denote
$$ e^{pq} \equiv e^{pq}(0) (p \neq q), e_*^{pq} \equiv e_*^{pq}(0), 
e_{**}^{pq} \equiv e_{**}^{pq}(0) (p \neq q), \/
p, q = 1, \cdots, l. $$ 
The operators $e^{pq}, e_*^{pq}, e_{**}^{pq} \;(p, q = 1, \cdots, l)$ form
the horizonal subalgebra $\frak {so}(2l)$ in $\widehat{ \frak {so}}(2l)$.
In particular, the operators $e_*^{pq}\; ( p, q = 1, \cdots, l)$ form
a subalgebra ${\frak {gl}}(l)$ in the horizontal ${\frak {so}}(2l)$.
We identify the Borel subalgebra ${\frak b}( {\frak {so}}(2l) )$
with the one generated by $e_{**}^{pq} \;(p \neq q), e_*^{pq}\; (p \leq q),
p, q = 1, \cdots, l.$

\begin{lemma}
  The action of ${\frak {gl}} (l)$ generated by
 $ e_*^{pq} \; ( p, q = 1, \cdots, l)$ 
  and that of $\hgl$ generated by 
  $ E_{ij}\; ( i,j \in \Z)$ on $\Fl$ commute with each other.
  \label{lem_basic}
\end{lemma}        
\begin{demo}{Proof}
One can prove $[e_*^{pq}, E_{ij}] = 0 $ by direct computation starting from
the definition of $e_*^{pq}$ and $E_{ij}$.
We will prove here the theorem in a more conceptual way by invoking the 
Wick's theorem. We will do it in the case
$ \underline{\Z} = \hz$. The case $ \underline{\Z} = \Z$
can be treated similarly by some slight modification.

The statement 1) is equivalent to 
\begin{equation}
  \left[
    \sum^l_{k=1} : \psi^{+,k} (z) \psi^{-,k} (w) :,
    \int : \psi^{+,p} (u) \psi^{-,q} (u):\, du
  \right] = 0, \quad p, q = 1, \dots l.
  \label{eq_zero}
\end{equation}
In order to prove (\ref{eq_zero}) we calculate some operator
product expansions (OPE) as follows.  Since
$$
\psi^{+,m} (z) \psi^{-,n} (w) \sim \frac{\delta_{m,n}}{z-w}, \quad
\psi^{-,m} (z) \psi^{+,n} (w) \sim \frac{\delta_{m,n}}{z-w}.
$$
we obtain by the Wick theorem
\begin{eqnarray*}
   \left(
    \sum^l_{k=1} : \psi^{+,k} (z) \psi^{-,k} (w) :
   \right)
   \left(
    : \psi^{+,p} (u) \psi^{-,q} (u):
  \right)                     
                    \\
   \quad  \quad \quad  \quad \sim 
  \frac{:\psi^{+,p} (z) \psi^{-,q} (u):}{w - u} +
  \frac{:\psi^{-,q} (w) \psi^{+,p} (u):}{z - u}.
\end{eqnarray*}
But for local fields $a (z)\/$ and $\psi^{-}(z)\/$ with OPE $a (z) b(u)
\sim \sum_j c_j (z) / (z - u)^j\/$ we have 
$\left[a (z), \int b(u)\, d u \right] = - c_1 (z)\/$.  
Hence the left-hand side of (\ref{eq_zero}) is equal to 
$:\psi^{+,p} (z) \psi^{-,q} (w): + \psi^{-,q} (w) \psi^{+,p} (z): = 0\/$.
\end{demo}
In the remaining part of this section, we will divide into
three cases: untwisted with $\underline{\Z} = \hz $,
untwisted with $\underline{\Z} = \Z $, and twisted cases.
\subsection{Untwisted case $\underline{\Z} = \hz$: dual pair $\Fpair$}
 Let
 \begin{eqnarray}
   \begin{array}{rcl}
     \lefteqn{
      \sum_{i,j \in \Z} ( E_{ij} - E_{1-j,1-i} ) z^{i- 1} w^{-j} =
             }          \\
    & & \sum^l_{k=1} (: \psi^{+,k} (z) \psi^{-,k} (w) : 
    - : \psi^{+,k} (w) \psi^{-,k} (z) :).
   \label{eq_gener}
   \end{array}
 \end{eqnarray}
We have the following lemma.
\begin{lemma}
    The action of the horizontal subalgebra ${\frak {so}}(2l)$
 and that of $\dinf$ generated by 
  $ E_{ij} - E_{1-j, 1-i} \;( i,j \in \Z)$ on $\Fl$ commute with each other.
 \label{lem_who}
\end{lemma}
\begin{demo}{Proof}
A similar argument to the proof of (\ref{eq_zero}) shows that 
\begin{eqnarray*}
  \left[
    \sum^l_{k=1} (: \psi^{+,k} (z) \psi^{-,k} (w) : 
   - : \psi^{+,k} (w) \psi^{-,k} (z) :),
   \right. \nonumber     \\
  \left.  \int : \psi^{+,p} (u) \psi^{+,q} (u):\, du
  \right] & =& 0, \\
  \left[
    \sum^l_{k=1} (: \psi^{+,k} (z) \psi^{-,k} (w) :   
    - : \psi^{+,k} (w) \psi^{-,k} (z) :), \right.  \nonumber     \\
  \left.  \int : \psi^{-,p} (u) \psi^{-,q} (u):\, du
  \right] & = & 0
  \label{eq_orth1}
 \end{eqnarray*}
 where $p, q = 1, \dots l, p \neq q.$
 The lemma now follows from the formula (\ref{eq_gener}).
\end{demo}
\begin{remark}
  We may introduce a natural $\Z_{+}/2$-gradation on
 $\Fl$ by the eigenvalues of a degree operator $d$ on $\Fl$
 which satisfies 
 $$ \left[ d, \psi_{-i}^{\pm, j} \right] = i \psi_{-i}^{\pm, j},
    \quad d \vac = 0. $$
 Each weight space of $\Fl$ is finite dimensional. Clearly
 an element in the horizontal subalgebra ${\frak {so}}(2l)$ has weight $0$.
 Thus as a representation of ${\frak {so}}(2l)$, each weight space
 of $\Fl$ is preserved by the action ${\frak {so}}(2l)$.
 So as a representation of ${\frak {so}}(2l)$, $\Fl$ is decomposed
 into a direct sum of finite dimensional irreducible representations. 
 As a representation of ${\frak {so}}(2l)$, 
 $\Fl$ is isomorphic to $ \wedge^* ( V \bigotimes \C^{\Bbb N} )$,
 where $V \cong \C^{2l}$ carries a natural symmetric bilinear form
 and $ \C^{ \Bbb N}$ is the direct limit of $\C^N$ as $N$ tends to 
 $ + \infty$. The action of the Lie algebra
 ${\frak {so}}(2l)$ can be lifted to an action of $SO(2l)$
 and naturally extends to $O(2l)$. 
   \label{rem_grad}
\end{remark}
Note that $\tau \in O(2l) - SO(2l)$ as defined in (\ref{eq_tau}) 
commutes with the Fourier components of  
$\psi^{\pm, k} (z), k = 1, \dots, l-1$, and sends 
$\psi^{+,l} (z)$ (resp. $\psi^{-,l} (z)$)
to $\psi^{-,l} (z)$ (resp. $\psi^{+,l} (z)$).
It follows  that $\tau$ commutes with
the Fourier components of the generating function
$\sum^l_{k=1} (: \psi^{+,k} (z) \psi^{-,k} (w) :
 - : \psi^{+,k} (w) \psi^{-,k} (z) :)$
which span $\dinf$ on $\Fl$ by equation (\ref{eq_gener}). Since
$g$ and $SO(2l)$ generate $O(2l)$, the following lemma follows 
from Lemma \ref{lem_who}.
\begin{lemma}
  The action of $O(2l)$ commutes with the action of
 $\dinf$ on $\Fl$.
\end{lemma}
Now we need to quote a lemma from the classical invariant theory 
(cf. e.g. \cite{H1}).
\begin{lemma}
   Let $G$ be the orthogonal group $O(k)$, and let $U$ and $W$ 
 be $G$-modules formed by taking direct sums of the 
 natural module for $G$. Consider the resulting action on 
 the tensor product $S(U) \bigotimes \wedge (V)$
 of the symmetric tensor $S(U)$ and anti-symmetric tensor $\wedge(V)$. 
Then the algebra of $G$-invariants 
 is generated by the invariants of degree $2$.
   \label{lem_invar}
\end{lemma}

It is easy to check that all the invariants of degree $2$
in $\Fl$ are precisely vectors obtained by letting elements
of $\dinf$ acting on the vacuum vector $\vac$. 
So we have a dual pair $\Fpair$ in the sense
of Howe acting on $\Fl$ (also cf. \cite{KR2}
for an appropriate however straightforward adaption of classical 
dual pairs to infinite dimensional Fock representation cases).
It can be argued similarly that $\Fpairgl$ also form a dual pair on $\Fl$.

Denote
\begin{eqnarray}
\Xi^{+, m}_i & \equiv & \psi^{+,i}_{-{m}+ \hf} 
\cdots \psi^{+,i}_{-\frac32} \psi^{+,i}_{-\hf},               \\  
\Xi^{-, m}_i & \equiv & \psi^{-,i}_{-{m}+ \hf} 
\cdots \psi^{-,i}_{-\frac32} \psi^{-,i}_{-\hf},               \\
\Xi_i^{ det} & \equiv & \psi^{+,i+1}_{- \hf} \psi^{-,i+1}_{- \hf}
\psi^{+,i+2}_{- \hf} \psi^{-,i+2}_{- \hf}  \ldots 
\psi^{+,l}_{- \hf} \psi^{-,l}_{- \hf}.
\end{eqnarray}
We take the convention that $\Xi^{\pm, 0}_i = 1$.
Define a map 
$\Lambda^{ \frak{aa}}_{+}: \Sigma (A) \longrightarrow \hgl_0^*$: 
$$\lambda = (m_1, \cdots, m_l) 
 \longmapsto \Lambda^{ \frak{aa}}_{+} (\lambda) $$
to be
 \begin{equation}
  \Lambda^{ \frak{aa}}_{+} (\lambda) = \hL_{m_1} + \cdots + \hL_{m_l}.
   \label{map_plus}
 \end{equation}
 We now recall a duality theorem for the dual 
pair $\Fpairgl$ which was earlier proved in \cite{F1, F2} 
(also see \cite{FKRW}). For the sake of completeness,
we sketch a proof here by invoking the dual pair principle 
in the spirit of this paper. Similar arguments will be needed
for the proofs of other duality theorems in this paper.
\begin{theorem}
  \begin{enumerate}
  \item[1)] We have the following  $(GL(l), \hgl )$-module
 decomposition:
 \begin{eqnarray*}
   \Fl = \bigoplus_{\lambda \in \Sigma (A)} 
          V( {\frak {gl}}(l); \lambda) \otimes L\left(
                                \hgl; \Lambda^{ \frak{aa}}_{+}(\lambda), l
                               \right)
 \end{eqnarray*}
 where $V( GL(l); \lambda)$ is
 the irreducible $GL (l)$-module
 of highest weight $\lambda$, and 
 $L( \hgl; \Lambda^{ \frak{aa}}_{+}(\lambda), l )$ 
 is the irreducible highest weight $\hgl$-module 
 of highest weight $\Lambda^{ \frak{aa}}_{+}(\lambda)$ and central charge $l$.

 \item[2)] Given $\lambda = (m_1, \ldots, m_l ) \in \Sigma (A)$, 
 assume that $$ m_1 \geq \cdots m_i \geq m_{i+1} = \cdots  = m_{j-1} = 0
 > m_j \geq \cdots \geq m_l. $$
 Then the highest weight vector corresponding to the weight 
 $\lambda \in \Sigma (A)$ is 
 \begin{eqnarray}
  v_{\Lambda} = \Xi_1^{+, m_1 } \ldots \Xi_i^{+, m_i } 
  \Xi_j^{-, -m_j } \ldots \Xi_l^{-, -m_l } \vac .   
   \label{eq_hvect}
 \end{eqnarray}
 \end{enumerate}
  \label{th_Fpairgl}
 \end{theorem}
\begin{demo}{Proof}
   First we can easily check that the vector $v_{\Lambda}$ in
 (\ref{eq_hvect}) is indeed a highest weight vector for
 $GL(l)$ and $\hgl$ respectively. Another direct calculation
 shows that the highest weight of $v_{\Lambda}$ for $GL(l)$
 is $ (m_1, \ldots, m_l )$. 

 The highest weight of $v_{\Lambda}$ for $\hgl$ can be
 read off from a table below conveniently.
 Similar tables will be used throughout this paper
 so we make some general remarks and conventions here.
 In the first row $k$ decreases one by one from left
 to right while in the second row the corresponding weights when 
 acting on $E_{k,k}$ 
 are listed. Weights are the same within each box in the second row.
 For $k$ greater than the first entry (which is $m_1$ in this case) 
 or less than the last entry in the first row 
 (which is $m_l -1$ in this case), the weight on the corresponding
 $E_{k,k}$ is zero. {\bf We further assume that 
 $ m_1 > \dots > m_i > m_{ i+1} = \ldots = m_j = 0
 > m_{j +1} >  \ldots >  m_l$ for the sake of simplicity of presenting
 the following table as a general convention used
 throughout the paper}: 

 \vspace{0.1in}
 \begin{center}
 \begin{tabular}
 {|l|ll|l|ll|lll|l|lll|} \hline
 $k$      & $m_1$, & \ldots & . & $m_i$  & \ldots   
   & 0,    & \ldots, & $m_{j+1}-1$& . & $m_{l-1}$,& \ldots & $m_l -1$  \\ 
  \hline
 $ E_{kk}$&   1,   & \ldots & . &   i    & \ldots  
   & $j-l$,& \ldots, &  $j-l$     & . & $ -1$,    & \ldots  & $-1$    \\ 
  \hline

 \end{tabular}
 \end{center}
 \vspace{0.1in}
    
 The general case 
 $ m_1 \geq \dots \geq m_i > m_{ i+1} = \ldots = m_j = 0
 > m_{j +1} \geq \ldots \geq m_l$ 
 can be treated in a similar way and the highest weight for the
 highest weight with respect to $\hgl$ can be shown easily 
 to be also $\Lambda^{ \frak{aa}}_{+} (\lambda)$.

 Now the part 1) follows from the general principle of
 dual pairs since we have found in 2) highest weight vectors
 of all the irreducible representations of $GL(l)$.
\end{demo}
\begin{remark}
   The irreducible representations
 $L( \hgl; \Lambda^{ \frak{aa}}_{+}(\lambda), l )$ exhaust all
 irreducible unitary representations of $\hgl$ of central charge $l$
 as $\lambda$ ranges over $\Sigma(A)$.
\end{remark}    
We define a map 
$\Lambda^{ \frak{dd}}: \Sigma (D) \longrightarrow {\dinf}_0^* $
(see Section \ref{subsec_evenorth} for the definition $\Sigma (D)$)
by sending $ \lambda = (m_1, \cdots, \overline{ m}_l) \; (m_l > 0 )$ to 
$$  \Lambda^{ \frak{dd}} (\lambda) =
   (l -i ) \hL_0^d + ( l -i ) \hL_1^d + \sum_{k =1}^i \hL_{m_k}^d, $$
sending $(m_1, \cdots, m_j, 0, \ldots, 0 )\; (j <l )$ to 
$$  \Lambda^{ \frak{dd}} (\lambda) =
   (2l -i -j) \hL_0^d + (j-i ) \hL_1^d + \sum_{k =1}^i \hL_{m_k}^d, $$
and sending $ (m_1, \ldots, m_j, 0, \ldots, 0 ) \bigotimes {det} 
   \; ( j < l)$ to 
 $$  \Lambda^{\frak{dd}} (\lambda) =
  (j-i) \hL_0^d + (2l -i -j) \hL_1^d + \sum_{k =1}^i \hL_{m_k}^d, $$
if $m_1 \geq \ldots m_{i} \geq m_{i+1} = \ldots = m_{j} =1
   > m_{j+1} = \ldots = m_l = 0.$

\begin{theorem}
 \begin{enumerate}
  \item[1)] We have the following 
 $(O(2l), \dinf )$-module decomposition:
 \begin{eqnarray*}
   \Fl = \bigoplus_{\lambda \in \Sigma (D) } I_{\lambda}
       \equiv \bigoplus_{\lambda \in \Sigma (D) } 
          V(O(2l); \lambda) \otimes L \left(\dinf;
                                 \Lambda^{\frak{dd}} (\lambda), l
                               \right)
 \end{eqnarray*}
 where $V(O(2l); \lambda)$ is the irreducible $O(2l)$-module
 parametrized by $\lambda \in \Sigma (D)$ and 
 $L \left(\dinf; \Lambda^{\frak{dd}} (\lambda), l  \right)$
 is the irreducible 
 highest weight $\dinf$-module of highest weight
 $ \Lambda^{\frak{dd}} (\lambda)$ and central charge $l$.
  \item[2)] With respect to $ ( {\frak {so} } (2l), \dinf )$, 
  \begin{enumerate}
  \item[a)] 
   the isotypic subspace $ I_{\lambda} $
  is decomposed into a sum of two irreducible representations
  with highest weight vectors  
  \begin{equation}
  \Xi_1^{+, m_1 } \cdots \Xi_{l-1}^{+, m_{l-1} } \Xi_l^{+, m_l }  \vac 
        \label{eq_h1}
 \end{equation}
  and
 \begin{equation}
  \Xi_1^{+, m_1 } \cdots \Xi_{l-1}^{+, m_{l-1} } \Xi_l^{-, m_l } \vac
  \label{eq_h2}
 \end{equation}
 in the case $\lambda = (m_1, \ldots, \overline{ m}_l )
 \in \Sigma(D), m_l > 0.$ 
 The highest weight of (\ref{eq_h1}) for ${\frak {so} } (2l)$
 is $(m_1, \ldots, m_{l-1}, m_l)$
 while that of (\ref{eq_h2}) for ${\frak {so} } (2l)$
 is $(m_1, \ldots, m_{l-1}, - m_l)$;
 \item[b)] 
  the isotypic subspace $ I_{\lambda} $ is irreducible
  with highest weight vector
  \begin{equation}
   \Xi_1^{+, m_1 } \cdots \Xi_j^{+, m_j}  \vac 
        \label{eq_h3}
 \end{equation}
 in the case $ \lambda = (m_1, \cdots, m_j, 0, \ldots, 0),
 m_1 \geq \dots \geq m_j >m_{j+1} = \ldots = m_{l} = 0, j <l; $
  \item[c)] 
  the isotypic subspace $ I_{\lambda} $ is irreducible
  with highest weight vector
  \begin{eqnarray}
   \Xi_1^{+, m_1 } \cdots \Xi_j^{+, m_j } \Xi_j^{det} \vac 
  \label{eq_h4}
 \end{eqnarray}
 in the case $ \lambda = (m_1, \cdots, m_j, 0, \ldots, 0 ) \bigotimes {det},$ 
 $ m_1 \geq \dots \geq m_j > 0, m_{j+1} = \ldots = m_l = 0, j < l. $
   \end{enumerate}
 \end{enumerate}
  \label{th_Fpair}
 \end{theorem}
\begin{demo}{Proof}
 We prove first that the vectors given in part 2) above indeed are of
 highest weight. The vectors in 2a) and 2b) are of highest weight
 for $\hgl$ and so they are also highest weight vectors for
 the Lie subalgebra $\dinf$. One can also prove by a direct computation
 that the vector in 2c) is also a highest weight vector for
 $\dinf$ (note that it is not a highest weight vector for $\hgl$
 so we cannot use the argument for 2a) and 2b)).   
 On the other hand, these vectors are of highest weight 
 for $GL (l) \subset O(2l)$ by Theorem \ref{th_Fpairgl}.
 Therefore to prove they are of highest weight for ${\frak {so} }(2l)$
 it suffices to show that they are also annihilated by 
 $$ e^{pq}_{**} = \SUM :\psi^{+,p}_{-n} \psi^{+,q}_n:,
    \quad p,q = 1, \dots, l, \; p \neq q. $$
 This is so for the vector (\ref{eq_h1}) 
 since both $\psi^{+,p}_{-n}$ and $ \psi^{+,q}_n$
 anticommute with any $\psi^{+,k}_i$ in 
 (\ref{eq_h1}) and thus one can move (up to a sign) 
 the annihilator between $\psi^{+,p} _{-n}$ and $\psi^{+,q}_n$ to 
 the right to kill the vacuum vector $\vac$.
 Noting that $ \left( \psi^{+,p}_{-n} \right)^2 = 0$,
 we also easily check that each term in the sum
 $e^{pq}_{**} = \SUM :\psi^{+,p}_{-n} \psi^{+,q}_n: \; 
 ( p,q = 1, \dots, l-1 )$ annihilates the vector (\ref{eq_h2})
 and so does $e^{pq}_{**}$.
 The case of (\ref{eq_h3}) is proved in the same way as in the 
 case of (\ref{eq_h1}).
 The case of (\ref{eq_h4}) also follows by noting that
 $ ( :\psi^{+,p}_{-n} \psi^{+,q}_n: ) \psi^{+,i}_{-\hf}
   \psi^{-,i}_{-\hf} \vac = 0$ for
   $p \neq q, p, q, i =1, \ldots, l.$

 It is easy to check case by case
 that the corresponding highest weights of these vectors
 in 2) above with respect to the Cartan subalgebra
 $ {\frak h} ({\frak {so} }(2l))$ 
 of ${\frak {so} }(2l)$ are given as in the theorem. 

 On the other hand, in the case of
 (\ref{eq_h3}), with $m_1 \geq \ldots m_{i} \geq m_{i+1} = \ldots = m_{j} =1
   > m_{j+1} = \ldots = m_l = 0,$ 
 we can easily calculate the highest weight when acting 
 on $E_{i,i}$ as in the following table. 

 \vspace{0.1in}
 \begin{center}
 \begin{tabular}
 {|l|ll|ll|l|lll|l|} \hline
 $k$      & $m_1$, & \ldots    & $m_2$,  & \ldots   
   & \ldots    & $m_i$, & \ldots, &  2  &  1  \\ \hline
 $ E_{kk}$ &   1,  & \ldots    &   2,    & \ldots  
   & \ldots    &  $i$,  & \ldots, & $i$ & $j$    \\ \hline

 \end{tabular}
 \end{center}
 \vspace{0.1in}

 Then we easily read off the highest weight with respect to
 $\dinf$ from the above table (see subsection \ref{subsec_dinf}
 for the definition of $\hL_i^d$):
 $$ \Lambda^{\frak{dd}} (\lambda) =
  (2l -i -j) \hL_0^d + (j-i ) \hL_1^d + \sum_{a =1}^i \hL_{m_a}^d . $$
 The case of (\ref{eq_h1}) can be
 treated as a special case $ j =l$ by using the preceding argument.
 Therefore we read off the corresponding highest weight for $\dinf$ as
 \begin{equation}
  \Lambda^{\frak{dd}} (\lambda) =
  (l -i) \hL_0^d + (l -i ) \hL_1^d + \sum_{a =1}^i \hL_{m_a}^d. 
   \label{eq_weightdd}
 \end{equation}
 Case of (\ref{eq_h2}) is divided into two subcases:
 first if $m_1 \geq \ldots \geq m_{l-1} \geq 2 \geq -m_l > 0$, we have
 the following table:

 \vspace{0.1in}
 \begin{center}
 \begin{tabular}{|l|ll|ll|l|lll|lll|} \hline
 $k$      & $m_1$, & \ldots    & $m_2$, & \ldots   
 & \ldots  & $m_{l-1}$, & \ldots, &  1   & 0, & \ldots,& $ 1-m_l $ 
     \\ \hline
 $E_{kk}$ &   1,   & \ldots    &   2,   & \ldots  
 & \ldots  & $l-1$,     & \ldots,  &      & $-1$,& \ldots,& $-1$        
     \\ \hline
 \end{tabular}
 \end{center}
 \vspace{0.1in}

 Secondly, if $m_1 \geq \ldots \geq m_i 
  > m_{i+1} = \ldots = m_l =1 $, we have the following table:

 \vspace{0.1in}
 \begin{center}
 \begin{tabular}{|l|ll|ll|l|lll|l|l|} \hline
 $k$      & $m_1$, & \ldots    & $m_2$, & \ldots   
   & \ldots    & $m_i$, & \ldots, &  2  & 1     & 0   \\ \hline
 $E_{kk}$ &   1,   & \ldots    &   2,   & \ldots  
   & \ldots    & $i$,   & \ldots, & $i$ & $l-1$ & -1   \\ \hline
 \end{tabular}
 \end{center}
 \vspace{0.1in}

 In either subcase,
  we easily obtain the highest weight with respect to $\dinf$ from the
 above two tables in the uniform formula (\ref{eq_weightdd}).

 In the case of (\ref{eq_h4}), we have

 \vspace{0.1in}
 \begin{center} \begin{tabular}{|l|ll|ll|l|lll|l|l|} \hline
 $k$      & $m_1$, & \ldots    & $m_2$, & \ldots   
   & \ldots    & $m_i$, & \ldots, &  2  & 1   &   0     \\ \hline

 $E_{kk}$ &   1,   & \ldots    &   2,   & \ldots  
   & \ldots    & $i$,   & \ldots, & $i$ & $l$ & $j-l$   \\ \hline
 \end{tabular}
 \end{center}
 \vspace{0.1in}

 Then we read off the highest weight in this case:
 $$  \Lambda^{\frak{dd}} (\lambda) =
  (j-i) \hL_0^d + (2l -i -j) \hL_1^d + \sum_{a =1}^i \hL_{m_a}^d. $$

 Now the decomposition in part 1) follows from the general
 abstract nonsense of dual pair theory since 
 we have found highest weight vectors of 
 all the irreducible representations of $O(2l)$ in 2).
\end{demo}
\begin{remark}
   Irreducible representations
 $L( \dinf; \Lambda^{ \frak{dd}}_{+}(\lambda), l )$ exhaust all
 irreducible unitary representations of $\dinf$ of central charge $l$
 as $\lambda$ ranges over $\Sigma(D)$.
\end{remark}    
We have the following corollary.
\begin{corollary}
   The space of invariants of $O(2l)$ in the Fock space $\Fl$ 
 is the irreducible 
 $\dinf$-module $ L (\dinf; 2l \hL^d_0 )$ of central charge $l$.
\end{corollary}
\begin{remark}
  The Dynkin diagram of $\dinf$ admits an automorphism of order 2 denoted
 by $\sigma$. $\sigma$ induces naturally an automorphism of order 2
 of $\dinf$, which is denoted again by $\sigma$ by abuse od notation.
 $\sigma$ acts on the set of highest weights of $\dinf$
 by mapping $ \lambda = h^d_0  \hL^d_0 + h^d_1  \hL^d_1 
  + \sum_{i \geq 2} h^d_i  \hL^d_i$ to
 $\sigma(\lambda ) = h^d_1  \hL^d_0 + h^d_0  \hL^d_1 
  + \sum_{i \geq 2} h^d_i  \hL^d_i$. 
 In this way one can obtain an irreducible module of the semi-product
 $\sigma \propto \dinf$
 on $ L(\dinf; \lambda) \bigoplus L(\dinf; \sigma(\lambda ) )$
 if $\sigma(\lambda ) \neq \lambda$ and on $ L(\dinf; \lambda) $
 if $\sigma(\lambda ) = \lambda $.
 
 Then it follows from Theorem \ref{th_Fpair} that the isotypic subspace
 of $\Fl$ with respect to the joint action 
 $( SO(2l), \sigma \propto \dinf)$ is irreducible and 
 $( SO(2l), \sigma \propto \dinf)$ form a dual pair on $\Fl$. In particular
 the space of invariants of $\Fl$ under the action fo $SO(2l)$ is 
 isomorphic to $L(\dinf; 2l \hL^d_0 ) \bigoplus L(\dinf; 2l \hL^d_1 )$.
   \label{rem_semi}
\end{remark}
\subsection{Untwisted case $\underline{\Z} = \Z$: dual pair $\Fpairdbtwo$}
  It follows from formula (\ref{eq_genelin}) that
 \begin{eqnarray*}
   \begin{array}{rcl}
    \lefteqn{
      \sum_{i,j \in \Z} ( E_{i,j} -  E_{-j,-i}) z^i w^{-j} 
            }     \\
     & = & \sum_{k =1}^l \left( : \psi^{+,k} (z) \psi^{-,k} (w):
                          - : \psi^{+,k} ( w) \psi^{-,k} ( z):
                  \right).
   \end{array}
 \end{eqnarray*}
\begin{lemma}
    The action of the horizontal subalgebra ${\frak {so}}(2l)$
  and that of $\binftwo$ generated by 
  $ E_{ij} - E_{-j, -i} \; ( i,j \in \Z )$ on $\Fl$ commute with each other.
 \label{lem_abel}
\end{lemma}
\begin{demo}{Proof}
We have shown in Lemma \ref{lem_basic}
that $E_{ij}, i,j \in \Z$ commutes with
the Lie subalgebra ${\frak {gl}}(l)$ of ${\frak {so}}(2l)$.
So it remains to prove that $E_{ij}$ commutes with
$e^{pq}_{**}, e^{pq}, p, q = 1, \ldots , l$.

\begin{eqnarray*}
 [E_{ij}, e^{pq}_{**}]
 & = & \left[ \sum_{k =1}^l :\psi^{+,k}_{-i} \psi^{-,k}_j :, 
 \sum_{ m \in \Z} :\psi_{-m}^{+,p} \psi_m^{+,q}: \right]          \\
 & = & \left[ \sum_{k =1}^l \psi^{+,k}_{-i} \psi^{-,k}_j , 
 \sum_{ m \in \Z} \psi_{-m}^{+,p} \psi_m^{+,q} \right]            \\
 & = & \left[ \psi^{+,p}_{-i} \psi^{-,p}_j,  \psi_{-j}^{+,p} \psi_j^{+,q} ]
 +  [ \psi^{+,q}_{-i} \psi^{-,q}_j, \psi_j^{+,p} \psi_{-j}^{+,q} \right]  \\
 & = & \psi^{+,p}_{-i} \psi_j^{+,q} - \psi^{+,q}_{-i} \psi_j^{+,p}.
\end{eqnarray*}
It follows immediately that
  $ [E_{ij}-E_{-j,-i}, e^{pq}_{**}] = 0. $
Similarly we can prove that
 $ [E_{ij}-E_{-j,-i}, e^{pq}] = 0. $
\end{demo}
\begin{remark}
  As a representation of ${\frak {so}}(2l)$, the Fock space
 $\Fl$ is isomorphic to $\wedge (\C^l )
 \wedge (\C^{2l} \bigotimes \C^{\Bbb N})$, where $\wedge (\C^l )$
 is the sum of two half-spin representations and 
 ${\frak {so}}(2l)$ acts on $\C^{2l} \bigotimes \C^{\Bbb N}$
 naturaly by the left action on $\C^{2l}$. The action of
 ${\frak {so}}(2l)$ can be lifted to $Spin(2l)$
 which extends naturally to $Pin (2l)$. It follows that any irreducible
 representation of $Spin(2l)$ appearing in $\Fl$ cannot
 factor to $SO(2l)$.  A similar
 argument to the classical dual pair case \cite{H2}
 shows that $Pin (2l)$ and $\binf$ form a dual pair on $\Fl$.
\end{remark}
Denote
\begin{eqnarray}
\Sigma^{+, m}_i & \equiv & \psi^{+,i}_{-{m}} 
\cdots \psi^{+,i}_{-2} \psi^{+,i}_{-1},         \nonumber       \\  
\Sigma^{-, m}_i & \equiv & \psi^{-,i}_{-m} 
\cdots \psi^{-,i}_{-1} \psi^{-,i}_{0}, \quad m \geq 0    \nonumber          
\end{eqnarray}
with the convention $\Sigma^{+, 0}_i = 1$.
We define a map $\Lambda^{\frak{db} }$
from $\Sigma(Pin) $ to ${\binf}_0^*$ by sending 
 $$ \lambda = (m_1, \ldots, \overline{ m}_l) , \quad
 m_1 \geq m_2 \geq \ldots \geq m_l \geq 0$$
to
 $$ \Lambda^{\frak{db}}(\lambda) =
     (2l - 2j) \hL_0^b + \sum_{k =1}^j \hL_{m_k}^b $$
if $m_1 \geq \ldots m_j > m_{j +1} = \ldots = m_l = 0.$

\begin{theorem}
 \begin{enumerate}
  \item[1)] We have the following 
  $\Fpairdbtwo$-module decomposition:
 \begin{eqnarray*}
  \Fl =  \bigoplus_{\lambda \in \Sigma (Pin) } I_{\lambda}
      \equiv \bigoplus_{\lambda \in \Sigma (Pin) } 
          V(Pin(2l); \lambda) \otimes L \left(\binftwo;
                                 \Lambda^{\frak{db}}(\lambda), l
                               \right)
 \end{eqnarray*}
 where $V(Pin(2l); \lambda)$ is the irreducible $Pin(2l)$-module
 parametrized by $\lambda \in \Sigma (Pin)$, and
 $L \left(\binftwo; \Lambda^{\frak{db}} (\lambda), l \right)$ 
 is the irreducible highest weight $\binftwo$-module of highest weight
 $ \Lambda^{\frak{db}} (\lambda)$ and central charge $l$.
 \item[2)] With respect to $( Spin(2l), \binftwo)$, 
  the isotypic subspace $I_{\lambda}$ is decomposed into
  a sum of two irreducible representations with highest weight vectors
  \begin{equation}
   \Sigma_1^{+, m_1 } \ldots \Sigma_{l-1}^{+, m_{l-1} }
  \Sigma_l^{+, m_l }  \vac  
        \label{eq_hig1}
 \end{equation}
  and
 \begin{equation}
   \Sigma_1^{+, m_1 } \ldots \Sigma_{l-1}^{+, m_{l-1} }
  \Sigma_l^{-, m_l } \vac 
  \label{eq_hig2}
 \end{equation}
  where $\lambda = \hf |{\bf 1_l }| + (m_1, \ldots, m_l) \in \Sigma (Pin),$
  $m_1 \geq \dots \geq m_{l-1} \geq m_l  \geq 0$.
  The highest weight of (\ref{eq_hig1}) for $Spin (2l)$ is 
  $ \hf {\bf 1}_l + (m_1, \ldots, m_l)$ while that of 
  (\ref{eq_hig2}) for $Spin (2l)$ is
  $ \hf {\bf \bar{1}}_l + (m_1, \ldots, m_{l-1}, -m_l)$.
 \end{enumerate}
 \label{th_pinpair}
\end{theorem}
\begin{demo}{Sketch of a proof}
   Proof is similar to that of Theorem \ref{th_Fpair}.
  Below we caculate the highest weights of the vectors (\ref{eq_hig1})
 and (\ref{eq_hig2}) for $\binftwo$. We list the weight on $E_{k,k}$
 for the vector (\ref{eq_hig1}) in the following table:
 
 \vspace{0.1in}
 \begin{center} \begin{tabular}{|l|ll|ll|l|lll|} \hline
 $k$      & $m_1$, & \ldots    & $m_2$,  & \ldots   
   & \ldots    & $m_j$, & \ldots, &  1   \\ \hline
 $ E_{kk}$ &   1,  & \ldots    &   2,    & \ldots  
   & \ldots    &  $j$,  & \ldots, & $j$   \\ \hline
 \end{tabular}
 \end{center}
 \vspace{0.1in}

 if $m_1 \geq \cdots \geq m_j > m_{j +1} = \cdots = m_l = 0.$
 Then we can easily read off the highest weight of the vector
 (\ref{eq_hig1}) for $\binftwo$ as follows:
 $$\Lambda^{\frak{db}} (\lambda) =    
       (2l - 2j) \hL_0^b + \sum_{k =1}^j \hL_{m_k}^b . $$

 Case of (\ref{eq_hig2}) is divided into two subcases:
 first if $m_1 \geq \ldots \geq m_j > m_{j +1} = \ldots = m_l = 0, j <l$,
 we have the following table:

 \vspace{0.1in}
 \begin{center} \begin{tabular}{|l|ll|ll|l|lll|l|} \hline
 $k$      & $m_1$, & \ldots    & $m_2$,  & \ldots   
   & \ldots    & $m_j$, & \ldots, &  1  & 0  \\ \hline
 $ E_{kk}$ &   1,  & \ldots    &   2,    & \ldots  
   & \ldots    &  $j$,  & \ldots, & $j$ & $-1 $   \\ \hline
 \end{tabular}
 \end{center}
 \vspace{0.1in}

 From this, we can read off the highest weight of the vector
 (\ref{eq_hig2}) for $\binftwo$ as follows:
 $$\Lambda^{\frak{db}} (\lambda) =    
       (2l - 2j) \hL_0^b + \sum_{k =1}^j \hL_{m_k}^b. $$
 Secondly, if $m_l \geq 1$, we have the following table:
 
 \vspace{0.1in}
 \begin{center}
 \begin{tabular}{|l|ll|ll|l|lll|lll|} \hline
 $k$       & $m_1$, & \ldots    & $m_2$,  & \ldots   
   & \ldots & $m_{l-1}$, & \ldots, &  1    & 0,   & \ldots & $-m_l$  \\ \hline
 $ E_{kk}$ &   1,   & \ldots    &   2,    & \ldots  
   & \ldots &  $l-1$,    & \ldots, & $l-1$ & $-1$,& \ldots & $-1$    \\ \hline
 \end{tabular}
 \end{center}
 \vspace{0.1in}

 From this, we can see that the highest weight of the vector
 (\ref{eq_hig2}) for $\binftwo$ to be
 $ \Lambda^{\frak{db}} (\lambda) = \sum_{k =1}^l \hL_{m_k}^b. $
\end{demo}
\begin{remark}
   The irreducible representations
 $L( \binftwo; \Lambda^{ \frak{db}}_{+}(\lambda), l )$ exhaust 
 irreducible unitary representations of $\binftwo$ of central charge $l$
 as $\lambda$ ranges over $\Sigma(Pin)$.
\end{remark}    
\subsection{Twisted case: dual pairs $\Fpaircc$ and $\Fpairdb$}
   It is well known \cite{FF} that the Fourier components of 
the following ``twisted'' generating functions 
\begin{eqnarray}
  \tilde{e}^{pq} (z) & \equiv &
           \sum_{n \in \Z} \tilde{e}^{pq} (n) z^{ -n -1 +2 \epsilon}
    =  : \psi^{-,p}(z) \psi^{-,q} ( -z): 
    ,  \nonumber   \\ 
  \tilde{e}_{**}^{pq} (z) & \equiv &
           \sum_{n \in \Z} \tilde{e}_{**}^{pq}(n) z^{ -n -1 +2 \epsilon}
    =  : \psi^{+,p}(z) \psi^{+,q} ( -z):,  \label{eq_twistfcn} \\
  e_*^{pq} (z) & \equiv & \sum_{n \in \Z} e_*^{pq} (n) z^{ -n -1 +2 \epsilon}
    =  : \psi^{+,p}(z) \psi^{-,q} (z): +  \delta_{p,q} \epsilon z^{-1},
   \quad p,q = 1, \dots, l   \nonumber    
\end{eqnarray}
span a representation of the twisted affine algebra
$ {\frak gl}^{(2)} ( 2l)$ 
of type $ A^{(2)}_{2l-1}$ with central charge $ 1$. 
Denote
\begin{eqnarray}
 \tilde{e}^{pq} \equiv \tilde{e}^{pq}(0), e_*^{pq} \equiv e_*^{pq}(0), 
  \tilde{e}_{**}^{pq} \equiv \tilde{e}_{**}^{pq}(0), \;
 p, q = 1, \cdots, l.
  \label{eq_horiz}
\end{eqnarray}
Now we divide into two cases according to $\underline{\Z} =\hz$ or $\Z$. 

\vspace{0.15in}
\noindent{\bf I. Case ${\bf \underline{\Z} =\hz}$: dual pair ${\bf \Fpaircc}$}
\vspace{0.1in}

It is easy to check that 
$ \tilde{e}^{pq} = \tilde{e}^{qp}$ and 
$ \tilde{e}_{**}^{pq} = \tilde{e}_{**}^{qp}$.
The horizontal subalgebra of $ {\frak gl}^{(2)} ( 2l)$ spanned by
the operators 
$\tilde{e}^{pq}, e_*^{pq}, \tilde{e}_{**}^{pq}, \; (p, q = 1, \cdots l) $
is isomorphic to Lie algebra $\frak {sp}(2l)$.
In particular, the operators $e_*^{pq}\; ( p, q = 1, \cdots l)$ form
a subalgebra ${\frak {gl}}(l)$ in the horizontal ${\frak {sp}}(2l)$.
We identify the Borel subalgebra ${\frak b}( {\frak {sp}}(2l) )$
with the one generated by $e_*^{pq}\; (p \leq q), \tilde{e}_{**}^{pq}, \;
p, q = 1, \ldots, l.$ Let
 \begin{eqnarray*}
   \begin{array}{rcl}
    \lefteqn{
      \sum_{i,j \in \Z} ( E_{i,j} -  (-1)^{i +j} E_{1-j, 1-i}) z^{i-1} w^{-j} 
            }     \\
     & = & \sum_{k =1}^l \left( : \psi^{+,k} (z) \psi^{-,k} (w):
                          + : \psi^{+,k} ( -w) \psi^{-,k} ( -z):
                  \right).
   \end{array}
 \end{eqnarray*}
The following lemma is straightforward.
\begin{lemma}
    The action of the horizontal subalgebra ${\frak {sp}}(2l)$
  and that of $\cinf$ generated by 
  $ E_{ij} - (-1)^{i +j} E_{1-j, 1-i} \; ( i,j \in \Z )$
 on $\Fl$ commute with each other.
\end{lemma}
\begin{remark}
   The action of ${\frak {sp}}(2l)$ can be integrated to an action
  of $Sp (2l)$ on $\Fl$. $Sp (2l)$ and $\cinf$ form a dual
  pair on $\Fl$.
\end{remark}
We define a map $\Lambda^{\frak{cc} }$ from $\Sigma(C)$ to ${\cinf}_0^*$ by
sending $(m_1, \ldots, m_l)$ to
$$ \Lambda^{\frak{cc} } (\lambda) =
       (l -j ) {}^c \hL_0 + \sum_{k =1}^j {}^c \hL_{m_j} $$
where $ m_1 \geq \ldots \geq m_j > m_{j +1} = \ldots = m_l = 0.$
\begin{theorem}
 \begin{enumerate}
  \item[1)] We have the following 
  $\Fpaircc$-module decomposition:
 \begin{eqnarray*}
  \Fl =  \bigoplus_{\lambda \in \Sigma (C) } I_{\lambda}
      \equiv \bigoplus_{\lambda \in \Sigma (C) } 
          V(Sp( 2l); \lambda) \otimes L \left(\cinf;
                                 \Lambda^{ \frak{cc} }(\lambda), l
                               \right)
 \end{eqnarray*}
 where $V(Sp(2l); \lambda)$ is the irreducible $Sp(2l)$-module
 parametrized by $\lambda \in \Sigma (Sp)$ and
 $L \left(\cinf; \Lambda^{\frak{cc} } (\lambda), l \right)$ 
 is the irreducible highest weight $\cinf$-module of highest weight
 $ \Lambda^{\frak{cc}} (\lambda)$ and central charge $l$.
 \item[2)] With respect to $ (Sp (2l), \cinf )$, the highest weight 
 vector in the isotypic subspace $I_{\lambda}$ is:
  \begin{eqnarray}
   \Xi_1^{+, m_1 } \ldots \Xi_l^{+, m_l }  \vac.
     \label{eq_vecsimple}
  \end{eqnarray}
 \end{enumerate}
\end{theorem}
\begin{demo}{Proof}
   Proof is similar to that of Theorem~\ref{th_Fpair}. Below we calculate
 the highest weight of the vector (\ref{eq_vecsimple}) for $\cinf$
 from the following table:

 \vspace{0.1in}
 \begin{center} \begin{tabular}{|l|ll|ll|l|lll|} \hline
 $k$      & $m_1$, & \ldots    & $m_2$,  & \ldots   
   & \ldots    & $m_j$, & \ldots, &  1   \\ \hline
 $ E_{kk}$ &   1,  & \ldots    &   2,    & \ldots  
   & \ldots    &  $j$,  & \ldots, & $j$   \\ \hline
 \end{tabular}
 \end{center}
 \vspace{0.1in}

 From this we easily see that the highest weight of the vector 
 (\ref{eq_vecsimple}) for $\cinf$ is 
 equal to $\Lambda^{ \frak{cc} }(\lambda)$.
\end{demo}
\begin{remark}
   The irreducible representations
 $L( \cinf; \Lambda^{ \frak{cc}}_{+}(\lambda), l )$ exhaust 
 irreducible unitary representations of $\cinf$ of central charge $l$
 as $\lambda$ ranges over $\Sigma(C)$.
\end{remark}    
\vspace{0.15in}
\noindent{\bf II. Case ${\bf \underline{\Z} =\Z}$: dual pair ${\bf \Fpairdb}$}
\vspace{0.1in}

    In this case it is easy to check that
$ e^{pq} = - e^{qp}$ and $ e_{**}^{pq} = - e_{**}^{qp}$.
The horizontal subalgebra of $ {\frak gl}^{(2)} ( 2l)$ spanned by
the operators $e^{pq}, e_*^{pq}, e_{**}^{pq}, \; (p, q = 1, \cdots, l)$
is isomorphic to Lie algebra $\frak {so}(2l)$.
In particular, the operators $e_*^{pq}\; ( p, q = 1, \cdots l)$ form
a subalgebra ${\frak {gl}}(l)$ in the horizontal ${\frak {sp}}(2l)$.
We identify the Borel subalgebra ${\frak b}( {\frak {so}}(2l) )$
with the one generated by $e_*^{pq}\; (p \leq q), e_{**}^{pq}, \;
p, q = 1, \cdots, l.$ Let
 \begin{eqnarray*}
   \begin{array}{rcl}
    \lefteqn{
      \sum_{i,j \in \Z} ( E_{i,j} -  (-1)^{i +j} E_{-j, -i}) z^i w^{-j} 
            }     \\
     & = & \sum_{k =1}^l \left( : \psi^{+,k} (z) \psi^{-,k} (w):
                          - : \psi^{+,k} ( -w) \psi^{-,k} ( -z):
                  \right).
   \end{array}
 \end{eqnarray*}
\begin{lemma}
    The action of the horizontal subalgebra ${\frak {so}}(2l)$
  and that of $\binf$ generated by 
  $ E_{ij} - (-1)^{i +j} E_{-j, -i} \; ( i,j \in \Z )$
 on $\Fl$ commute with each other.
\end{lemma}
Proof is similar to that of Lemma \ref{lem_abel}.
\begin{remark}
  The action of ${\frak {so}}(2l)$ can be integrated to an action
 of $Spin (2l)$ and extends naturally to $Pin (2l)$ on $\Fl$. 
 Furthermore, $Pin (2)$ and $ \binf$ form a dual pair on $\Fl$.
 We can also obtain a duality theorem for the dual pair
 $\Fpairdb$. Indeed the theorem can be
 stated in the same way as Theorem~\ref{th_pinpair}
 by replacing $\binftwo$ in Theorem~\ref{th_pinpair}
 by $\binf$. We omit the detail here.
 The dual pair $\Fpairdb$ is related to the dual pair $\Fpairdbtwo$
 by choosing a different symmetric bilinear form in defining the
 action on $\Fl$ of an infinite dimensional Lie subalgebra
 of $\hgl$ of $B$ type (cf. Section \ref{subsec_binf}).
\end{remark}
\section{Duality in the Fock space $\Flpmhalf$}
 \label{sec_Fpairbd} 
 Introduce a neutral fermionic field 
$\phi (z) = \sum_{n \in \underline{\Z} } \phi_n z^{-n -\hf + \epsilon}$ 
which satisfies the following commutation relations:
$$[ \phi_m , \phi_n ]_{+} = \delta_{m, -n},
     \quad m,n \in \underline{\Z}.$$
Denote by $\Fhalf$ the Fock space of $\phi (z)$
generated by a vacuum vector $\vac$, which is annihilated by 
$\phi_n, n \in \underline{\Z}_{+}$.
Denote by $\Flhalf$ the tensor product of $\Fhalf$ 
and the Fock space $\Fl$ of $l$ pairs of fermionic fields 
$\psi^{\pm, k} (z) (k =1, \dots, l)$.

Denote by 
\begin{eqnarray}
 e^{p} (z) \equiv \sum_{n \in \Z}
 e^{p}(n) z^{-n -1 + 2 \epsilon} & = & : \psi^{-,p}(z) \phi (z):, 
             \nonumber       \\
 e_*^{p} (z) \equiv \sum_{n \in \Z}
 e_*^{p}(n) z^{-n -1 + 2 \epsilon} & = & : \psi^{+,p}(z) \phi (z):, 
    \quad p = 1, \dots, l.         \label{eq_oddorth}
\end{eqnarray}

Then the Fourier components of $e^p (z), e^p_* (z)$, together with
generating functions (\ref{eq_orthaffine}),
(\ref{eq_orth}) and (\ref{eq_orthaffine})
 $$ e^{pq}(n) (p \neq q),\; e_{**}^{pq}(n) (p \neq q),\; e_*^{pq}(n),\; 
e^{p} (n), e_*^{p}(n) \; \; ( n \in \Z, \; p, q = 1, \dots, l) $$
generate an affine algebra $\widehat{{\frak {so}}}(2l+1)$ of level $1$ 
\cite{F, KP}. Put
\begin{eqnarray*}
     e_{**}^{pq} \equiv e_{**}^{pq}(0) (p \neq q),
 & e_*^{p} \equiv e_*^{p}(0),& e_*^{pq} \equiv e_*^{pq}(0), \\
     e^{pq} \equiv e^{pq}(0) (p \neq q),
 & e^{p} \equiv e^{p} (0), &p, q = 1, \dots, l.
\end{eqnarray*}
Then 
$ e^{pq} (p \neq q), e_{**}^{pq}(p \neq q), e_*^{pq}, 
e^{p}, e_*^{p} \;( p, q = 1, \dots, l) $ 
generate the horizontal subalgebra 
${\frak {so}}(2l+1)$ of $\widehat{{\frak {so}} }(2l+1)$.
In particular, we identify the Borel subalgebra
${\frak b}( {\frak {so}}(2l +1) )$
with the one generated by $e_{**}^{pq} (p \neq q), e_*^{pq} (p \leq q), \;
e_*^p, p, q = 1, \cdots, l.$
\subsection{Untwisted case $ \underline{\Z} = \hf + \Z $: dual pair $\Fpairbd$}
\begin{lemma}
   \begin{enumerate}
 \item[1)] Putting
 \begin{eqnarray}
   \begin{array}{rcl}
     \lefteqn{
   \sum_{i,j \in \Z} ( E_{ij}
    - E_{1-j,1-i} ) z^{i- 1} w^{-j}
             }             \nonumber \\
   & & = \sum^l_{k=1} (: \psi^{+,k} (z) \psi^{-,k} (w) :
    - : \psi^{+,k} (w) \psi^{-,k} (z) :) 
    + :\phi (z) \phi (w):
   \end{array}
 \end{eqnarray}
 defines an action of $\dinf$ on $\Flhalf$ with central charge $l + \hf$.
 \item[2)] The action of the horizontal subalgebra ${\frak {so}}(2l +1)$
 commutes with that of $\dinf$ generated by 
 $ E_{ij} - E_{1-j, 1-i} \; (i,j \in \Z )$ on $\Flhalf$.
   \label{lem_oddabel}
 \end{enumerate}
\end{lemma}        
\begin{demo}{Proof}
 Part 1) can be proved by a direct calculation.

 In the same way as we proved Lemma \ref{lem_who} we see that 
 $$ e^{pq}, e_*^{pq} (p \neq q), 
 e_{**}^{pq}(p \neq q)\; (p, q = 1, \dots, l) $$
 commute with the action of $\dinf$. So it remains to
 check that 
 $ e^{p}, e_*^{p} (p, q = 1, \dots, l) $ also commute
 with the action of $\dinf$.
 This is equivalent to showing that
 \begin{eqnarray}
  \left[
    \sum^l_{k=1} (: \psi^{+,k} (z) \psi^{-,k} (w) :
  - : \psi^{+,k} (w) \psi^{-,k} (z) :) 
  + :\phi (z) \phi (w):,  \right.      \nonumber   \\
  \left.  \int : \psi^{+,p} (u) \phi (u):\, du
  \right] = 0,       \label{eq_comm1}               \\
  \left[
    \sum^l_{k=1} (: \psi^{+,k} (z) \psi^{-,k} (w) : 
  - : \psi^{+,k} (w) \psi^{-,k} (z) :)
  + :\phi (z) \phi (w): , \right.       \nonumber   \\
  \left.  \int : \psi^{-,p} (u) \phi (u):\, du
  \right] = 0              
  \label{eq_comm2}
 \end{eqnarray}
 where $ p = 1, \dots l.$

 We will only prove (\ref{eq_comm1}) since the proof 
 of (\ref{eq_comm2}) goes parallely.

 In order to prove (\ref{eq_comm1}) we calculate some operator
 product expansions (OPE). Since
 $$
  \psi^{+,m} (z) \psi^{-,n} (w) \sim \frac{\delta_{m,n}}{z-w}, \quad
  \phi (z) \phi (w) \sim \frac{\delta_{m,n}}{z-w}.
 $$
 we have by using the Wick theorem, 
 \begin{eqnarray*}
  \Big(
    \sum^l_{k=1} (: \psi^{+,k} (z) \psi^{-,k} (w) - : \psi^{+,k} (w) \psi^{-,k} (z) :) 
   : + :\phi (z) \phi (w): 
  \Big) \cdot             \nonumber        \\
  \left(
    : \psi^{+,p} (u) \phi (u):
  \right)                    \nonumber       \\
  \sim 
  \frac{:\psi^{+,p} (z) \phi (u):}{w - u}
  - \frac{:\psi^{+,p} (w) \phi (u):}{z - u}
  - \frac{:\phi (z) \psi^{+,p} (u):}{w - u} 
  + \frac{:\phi (w) \psi^{+,p} (u):}{z - u}.
 \end{eqnarray*}
 But for local fields $a (z)\/$ and $b (z)\/$ with OPE $a (z) b(u)
  \sim \sum_j c_j (z) / (z - u)^j\/$ we have 
 $ \left[a (z), \int b(u)\, d u \right] = - c_1 (z)\/$.  
 Hence the left-hand side of (\ref{eq_comm1}) is equal to 
 $$ :\psi^{+,p} (z) \phi (w):
   - :\psi^{+,p} (w) \phi (z):
   - :\phi (z) \psi^{+,p} (w): 
   + :\phi (w) \psi^{+,p} (z): 
 $$
 which is equal to $0$.
\end{demo}
\begin{remark}
 The Fock space $\Flhalf$ is isomorphic to
 $ \wedge^* \left( \C^{2l+1} \bigotimes \C^{\Bbb N} \right)$
 as a representation of ${\frak {so}}(2l + 1)$.
 We can lift the action of ${\frak {so}}(2l + 1)$ to $ SO (2l + 1) $
 and extends to $O (2l + 1) $ naturally. A particular element
 $g \in O(2l +1) - SO(2l +1)$ is $diag (1, \ldots, 1, -1)$.
 It commutes  with the Fourier components of
 $\psi^{\pm ,k} (z) \; (k = 1, \dots, l)$, and sends 
 $\phi (z)$ to $- \phi (z)$. 
\end{remark}
It follows that $g$ commutes with the Fourier components of
 $\sum^l_{k=1} (: \psi^{+,k} (z) \psi^{-,k} (w) :
   - : \psi^{+,k} (w) \psi^{-,k} (z) :) 
   + :\phi (z) \phi (w):$
 which span $\dinf$ on $\Flhalf$. Since $g$ ans $SO(2l +1)$
generate $O(2l+1)$, the following lemma follows from Lemma \ref{lem_oddabel}.
\begin{lemma}
  The action of $O(2l+1)$ commutes with the action of
 $\dinf$ on $\Flhalf$.
\end{lemma}
\begin{remark}
 We can easily check that indeed all the invariants of degree $2$
 in $\Flhalf$ are exactly the vectors obtained by letting elements
 of $\dinf$ acting on the vacuum vector $\vac$ of $\Flhalf$. So
 by results on classical invariant theory \cite{H1}
 we have a dual pair $\Fpairbd$ acting on $\Flhalf$. 
\end{remark}
We define a map $\Lambda^{ \frak{bd} } $ from $\Sigma (B)$ 
(see Section \ref{sec_classical} for the definition of $\Sigma(B)$)
to ${\dinf}_0^*$ by sending
$$\lambda = (m_1, m_2, \ldots, m_l)$$
 to
$$  \Lambda^{ \frak{bd} } (\lambda) =
     (2l+1 - i -j) \hL^d_0 + (j-i)\hL^d_1 + \sum_{k=1}^i \hL^d_{m_k}  $$
and sending
$$\lambda = (m_1, m_2, \ldots, m_l) \bigotimes det $$
 to
$$   \Lambda^{ \frak{bd} } (\lambda) =
    (j - i) \hL^d_0 + (2l +1-i -j)\hL^d_1 + \sum_{k =1}^i \hL^d_{m_k} $$
assuming that 
$$ m_1 \geq \ldots \geq m_i > m_{i +1} = \ldots = m_j = 1 >
m_{j+1} = \ldots = m_l = 0 .$$
\begin{theorem}
 \begin{enumerate}
  \item[1)] We have the following $\Fpairbd$-module
 decomposition:
 \begin{eqnarray*}
   \Flhalf = \bigoplus_{\lambda \in \Sigma(B) } I_{\lambda}
       \equiv \bigoplus_{\lambda \in \Sigma(B) } 
          V(O(2l+1); \lambda) \otimes L \left(\dinf;  
                                 \Lambda^\frak{bd} (\lambda), l +1/2
                               \right)
 \end{eqnarray*}
 where $V(O(2l+1); \lambda)$ is the irreducible $O(2l+1)$-module
 parametrized by $\lambda  \in \Sigma(B) $ and 
 $L \left(\dinf; \Lambda^\frak{bd} (\lambda), l+1/2 \right)$
 is the irreducible highest weight $\dinf$-module with 
 highest weight $ \Lambda^\frak{bd} (\lambda)$ and central charge $l +1/2$.
  \item[2)] With respect to $(${\frak {so}}(2l +1)$, \dinf) $, 
 the isotypic subspace $I_{\lambda}$ is decomposed into a sum of
two irreducible representations with  highest weights
  \begin{eqnarray}
    \Xi_1^{+,m_1} \ldots \Xi_l^{+, m_l} \vac
  \mbox{ for } \lambda =(m_1, m_2, \ldots, m_l),   \label{eq_high1}
  \end{eqnarray}
  and
  \begin{eqnarray}
   \Xi_1^{+,m_1 } \ldots \Xi_j^{+, m_j} \Xi_j^{det} \phi_{-\hf} \vac 
  \mbox{ for } \lambda =(m_1, m_2, \ldots, m_l) \bigotimes det.
     \label{eq_high2}
  \end{eqnarray}
 \end{enumerate}
  \label{th_Fpairbd}
 \end{theorem}
\begin{demo}{Sketch of a proof}
   Proof is again similar to that of Theorem \ref{th_Fpair}.
 We indicate below only the points which are different
 from that of Theorem \ref{th_Fpair}.
 The determination of the highest weight
 of the vector (\ref{eq_high1}) for $\dinf$ is the same as before.
 $E_{i,i} - E_{1-i, 1-i} \;( i \in \Bbb N)$ acts on the vector 
 by the operator $\phi_{-i + \hf}\phi_{i - \hf}$. Then it
 is easy to see that
 $ (E_{i,i} - E_{1-i, 1-i} ) \phi_{ -\hf} \vac
   = \delta_{i,1} \phi_{ -\hf} \vac. $ 
 From this and the computation for the case of (\ref{eq_high1}), we obtain
 the highest weight of the vector (\ref{eq_high2}).
\end{demo}
We immediately have the following corollary
\begin{corollary}
    The space of invariants of $O(2l +1)$ in $\Flhalf$ is isomorphic to
 the irreducible $\dinf$-module $L ( \dinf; (2l +1) \hL^d_0 )$
 of central charge $l + \hf$.
\end{corollary}
\begin{remark}
   The irreducible representations
 $L( \dinf; \Lambda^{ \frak{bd}}_{+}(\lambda), l +1/2 )$ exhaust all
 irreducible unitary representations of $\dinf$ of central charge $l +1/2$
 as $\lambda$ ranges over $\Sigma(B)$.
\end{remark}    
\begin{remark}
 $( SO(2l +1), \sigma \propto \dinf)$ form a dual pair on $\Fl$.
 In particular the space of invariants of $\Flhalf$ under the action of
 $SO(2l +1)$ is isomorphic to
 $L(\dinf; (2l+1) \hL^d_0 ) \bigoplus L(\dinf;  (2l+1) \hL^d_1 )$.
 cf. Remark \ref{rem_semi}.
\end{remark}    
\subsection{Untwisted case $ \underline{\Z} = \Z $: dual pair $\Fpairbbtwo$}
We first have the following lemma.
\begin{lemma}
   Putting 
 \begin{eqnarray*}
   \begin{array}{rcl}
    \lefteqn{
      \sum_{i,j \in \Z} ( E_{i,j} -  E_{-j,-i}) z^i w^{-j} 
            }     \\
     & = & \sum_{k =1}^l \left( : \psi^{+,k} (z) \psi^{-,k} (w):
                          - : \psi^{+,k} ( w) \psi^{-,k} ( z):
                  \right)
                  + :\phi (z) \phi (w):
   \end{array}
 \end{eqnarray*}
 defines a representation of $\binftwo$ on $\Flhalf$ of central
 charge $2l +1$.
\end{lemma}        
Proof is straightforward.
\begin{lemma}
  The action of the horizontal subalgebra ${\frak {so}}(2l +1)$
 and that of $\binftwo$ generated by 
 $ E_{ij} - E_{-j, -i} \; ( i,j \in \Z )$ on $\Flhalf$ commute with each other.
\end{lemma}        
Define a map $\Lambda^{\frak{bb} }$ from $\Sigma (PB)$ to 
${\binf}_0^*$ by sending
$\lambda = \hf {\bf 1}_l + (m_1, m_2, \ldots, m_l)$ to
$$  \Lambda^{\frak{bb} } (\lambda)
     = (2l +1 - 2j ) \hL^b_0 + \sum_{k=1}^j \hL^b_{m_k} $$
if $m_1 \geq \ldots \geq m_j > m_{j+1} = \ldots = m_l =0.$
                 
The Fock space $\Flhalf$ splits into a sum of
two subspaces $\Flhalf_e$ and $\Flhalf_o$, where
$\Flhalf_e$ consists all even vectors while $\Flhalf_o$
consists all odd vectors according to the $ \Z_2$ gradation
on the vector superspace $\Flhalf$. Each subspace is
clearly invariant under the action of ${\frak {so}}(2l +1)$.
\begin{remark}
  The action of ${\frak {so}}(2l +1)$ can be integrated to
 $Spin (2l +1)$ on $\Flhalf_e$ and $\Flhalf_o$ respectively.
 One has indeed dual pair $\Fpairbbtwo$ acting on each
 $\Flhalf_e$ and $\Flhalf_o$.
\end{remark}
\begin{theorem}
 \begin{enumerate}
  \item[1)] We have the following $\Fpairbbtwo$-module decomposition:
 \begin{eqnarray}
   \Flhalf_e = \bigoplus_{\lambda \in \Sigma(PB) } 
          V(Spin(2l+1); \lambda) \otimes L \left(\binftwo;  
                                 \Lambda^{\frak{bb}} (\lambda), l + \hf
                               \right) ,       \label{eq_decomp1}          \\
   \Flhalf_o = \bigoplus_{\lambda \in \Sigma(PB) } 
          V(Spin(2l+1); \lambda) \otimes L \left(\binftwo;  
                                 \Lambda^{\frak{bb}} (\lambda), l + \hf
                               \right)  \label{eq_decomp2}
 \end{eqnarray}
 where $V(Spin(2l+1); \lambda)$ is the irreducible $Spin(2l+1)$-module
 parametrized by 
 $\lambda = \hf {\bf 1}_l + (m_1, m_2, \ldots, m_l)$ and 
 $ L \left(\binftwo; \Lambda^{\frak{bb}} (\lambda), l + \hf \right)$ 
 is the irreducible 
 highest weight $\binftwo$-module with central charge $l + \hf$.
  \item[2)] With respect to $\Fpairbbtwo$,
 the highest weight vectors corresponding to
 the weight $\lambda \in \Sigma (PB)$
 in (\ref{eq_decomp1}) and (\ref{eq_decomp2}) are respectively
 \begin{eqnarray}
  \Sigma_1^{+,m_1} \ldots \Sigma_l^{+, m_l} \vac,       \\
  \Sigma_1^{+,m_1} \ldots \Sigma_l^{+, m_l} \phi_0 \vac. 
 \end{eqnarray}
 \end{enumerate}
   \label{th_Fpairbb}
 \end{theorem}
\begin{remark}
   The irreducible representations
 $L( \binftwo; \Lambda^{ \frak{bb}}_{+}(\lambda), l+1/2 )$ exhaust all
 irreducible unitary representations of $\binftwo$ of central charge $l +1/2$
 as $\lambda$ ranges over $\Sigma(PB)$.
\end{remark}    
\subsection{Twisted cases: $\Fpairospc$ and $\Fpairbb$}
\vspace{0.15in}
\noindent{\bf I. Case ${\bf  \underline{\Z} = \hz}$:
 dual pair ${\bf \Fpairospc}$}
\vspace{0.1in}

 Introduce a bosonic field
$ \chi (z) = \sum_{ n \in \hz} \chi_n z^{ -n - \hf}$ which 
satisfies the following commutation relations:
$$[ \chi_m , \chi_n ] = ( -1)^{m + \hf}\delta_{m, -n}, \quad m,n \in \hz .$$
Denote by $\Fminushalf$ the Fock space of $\chi (z)$ generated by
a vacuum vector which is annihilated by $\chi_n, n\in \hz_{+}$.
Let $\Flminushalf$ be the tensor product of the Fock space 
of $l$ pairs of fermionic fields $\psi^{\pm, k} (z) \; (k = 1, \ldots, l )$
and the Fock space $\Fminushalf$ of $\chi (z)$.

It is known \cite{FF} that
the Fourier components of the generating functions 
$\tilde{e}^{pq} (z), \tilde{e}_{**}^{pq} (z)$ and $ e_*^{pq} (z)$
defined in (\ref{eq_twistfcn}) together with 
the following generating functions
\begin{eqnarray*}
 \zeta (z) & \equiv & \sum_{n \in \Z}
            \zeta(n) \NN  = : \chi (z) \chi ( -z):,   \nonumber     \\
 \tilde{e}^{p} (z) & \equiv &\sum_{n \in \Z}
            e^{p}(n) \NN  =  : \psi^{-,p}(z) \chi ( -z):,   \nonumber     \\
 \tilde{e}_*^{p} (z) & \equiv & \sum_{n \in \Z}
           e_*^{p}(n) \NN  = : \psi^{+,p}(z) \chi (z):,   
                        \quad (p,q = 1, \dots, l)  \nonumber    
\end{eqnarray*}
span a representation of the affine algebra ${\frak {gl}}^{(2)}(1, 2l )$
of type $ A^{(2)}(0, 2l-1)$ with central charge $ 1$. 
Denote
$$ \tilde{e}^{p} \equiv \tilde{e}^{p}(0), 
 \tilde{e}_*^{p} \equiv \tilde{e}_*^{p}(0), 
 \tilde{e}^{pq} \equiv \tilde{e}^{pq}(0),
  e_*^{pq} \equiv e_*^{pq}(0), 
 \tilde{e}_{**}^{pq} \equiv \tilde{e}_{**}^{pq}(0), \;
p, q = 1, \cdots, l. $$ 
Easy to check that the horizonal subalgebra in ${\frak {gl}}^{(2)}(1, 2l )$
spanned by the operators 
$\tilde{e}^{p}, \tilde{e}_*^{p}, \tilde{e}^{pq}, e_*^{pq}, 
 \tilde{e}_{**}^{pq}, \;(p, q = 1, \cdots, l)$
is isomorphic to Lie superalgebra $\frak {osp}(1, 2l)$.
In particular, the operators $e_*^{pq}\; ( p, q = 1, \cdots, l)$ form
a subalgebra ${\frak {gl}}(l)$ in the horizontal ${\frak {osp}}(1, 2l)$.
We identify the Borel subalgebra ${\frak b}( {\frak {osp}}(1, 2l) )$
with the one generated by $e_{**}^{pq}, e_*^{pq}\; (p \leq q),\;
p, q = 1, \cdots, l.$
\begin{lemma}
   Putting 
 \begin{eqnarray*}
   \begin{array}{rcl}
    \lefteqn{
      \sum_{i,j \in \Z} ( E_{i,j} - ( -1)^{i+j} E_{1-j,1-i}) z^{i-1} w^{-j}
            }     \\
   & = & \sum_{k =1}^l \left( : \psi^{+,k} (z) \psi^{-,k} (w):
                          + : \psi^{+,k} ( -w) \psi^{-,k} ( -z):
                  \right)
                  + :\chi (z) \chi (-w):
  \end{array}
 \end{eqnarray*}
 defines a representation of $\cinf$ on $\Flminushalf$ of central
 charge $l - \hf$.
\end{lemma}        
Proof is straightforward.
\begin{lemma}
     The action of the horizontal subalgebra ${\frak {osp}}(1, 2l )$
 and that of $\cinf$ generated by 
 $ E_{ij} - ( -1)^{i+j} E_{1-j, 1-i}\;( i,j \in \Z)$ on $\Flminushalf$
 commute with each other.
\end{lemma}        
\begin{remark}
  The action of the horizontal subalgebra ${\frak {osp}}(1, 2l )$
 can be integrated to $Osp (1, 2l)$. $Osp (1, 2l)$
 and $ \cinf$ form a dual pair on $\Flminushalf$.
\end{remark} 
We define a map $\Lambda^{ \frak{ospc}}$ from
$\Sigma (Osp)$ to $ {\cinf}_0^*$ by sending
$\lambda = (m_1, \ldots, m_l )$ to
$$\Lambda^{ \frak{ospc}} (\lambda) =
   (l - \hf -j) \hL^c_0 + \sum_{k =1}^j \hL^c_{m_k}$$
if $m_1 \geq \ldots \geq m_j > m_{j+1} = \ldots = m_l =0.$
\begin{theorem}
 \begin{enumerate}
  \item[1)] We have the following 
  $\Fpairospc$-module decomposition:
 \begin{eqnarray*}
  \Flminushalf =  \bigoplus_{\lambda \in \Sigma (Osp) } I_{\lambda}
           \equiv \bigoplus_{\lambda \in \Sigma (Osp) } 
          V(Osp(1, 2l); \lambda) \otimes L \left(\cinf;
                                 \Lambda^{ \frak{ospc}}(\lambda), l - 1/2
                               \right)
 \end{eqnarray*}
 where $V(Osp(1, 2l); \lambda)$ is the irreducible $Osp(1, 2l)$-module
 parametrized by $\lambda \in \Sigma (Osp)$, and
 $L \left(\cinf; \Lambda^{ \frak{ospc}} (\lambda), l - 1/2 \right)$ 
 is the irreducible highest weight $\cinf$-module of highest weight
 $ \Lambda^{ \frak{ospc}} (\lambda)$ and central charge $l - 1/2$.
 \item[2)] With respect to $ \Fpairospc$, 
 the isotypic subspace $I_{\lambda}$ is irreducible with 
 highest weight vector
  \begin{eqnarray*}
   \Xi_1^{+, m_1 } \ldots \Xi_l^{+, m_l }  \vac.
  \end{eqnarray*}
 \end{enumerate}
\end{theorem}
\vspace{0.15in}
\noindent{\bf II. Case ${\bf \underline{\Z} = \Z}$: 
dual pair ${\bf \Fpairbb }$ }
\vspace{0.1in}

Introduce a fermionic field
 $\varphi (z) = \sum_{n \in \Z} \varphi_n z^{-n}$ which 
satisfies the following commutation relations:
\begin{equation}
   [ \varphi_m , \varphi_n ]_{+} = ( -1)^m \delta_{m, -n},
     \quad m,n \in \Z. 
  \label{eq_commrel}
\end{equation}
   In this case the Fock space $\Flhalf$ is the tensor product
of the Fock space of $l$ pairs of fermionic fields 
$\psi^{\pm, k} (z) \;( k = 1, \ldots, l )$
and the Fock space $\Fhalf$ of $\varphi (z)$
generated by a vacuum vector which is annihilated by
$ \varphi_m, m \in \Bbb N$.

The Fourier components of the following generating functions
$$ : \varphi (z) \varphi ( -z):, \quad : \psi^{-,p}(z) \varphi ( -z):,
   \quad : \psi^{+,p}(z) \varphi (z): $$
together with generating functions $ \tilde{e}^{pq} (z),
\tilde{e}_{**}^{pq} (z), e_*^{pq} (z)$ defined in (\ref{eq_twistfcn})
                 %
form an affine algebra of type $ A^{(2)}_{ 2l}$
on $\Flhalf$. The horizontal subalgebra of 
$ A^{(2)}_{ 2l}$ is isomorphic to ${\frak {so}}( 2l +1 )$.
         
One can prove the following lemmas by a direct computation.
\begin{lemma}
   Putting 
 \begin{eqnarray*}
   \begin{array}{rcl}
    \lefteqn{
     \sum_{i,j \in \Z} ( E_{i,j} - ( -1)^{i+j} E_{-j,-i}) z^i w^{-j}
            }     \\
   & = & \sum_{k =1}^l \left( : \psi^{+,k} (z) \psi^{-,k} (w):
                          + : \psi^{+,k} ( -w) \psi^{-,k} ( -z):
                  \right)
                  + :\varphi (z) \varphi (-w):
   \end{array}
 \end{eqnarray*}
 defines a representation of $\binf$ on $\Flhalf$ of central
 charge $l + 1/2$.
\end{lemma}        
\begin{lemma}
     The action of the horizontal subalgebra ${\frak {so}}( 2l +1 )$
 commutes with that of $\binf$ generated by 
 $ E_{i,j} - ( -1)^{i+j} E_{-j,-i}), i,j \in \Z$ on $\Flhalf$
 commute with each other.
\end{lemma}        
We can also show that the action of ${\frak {so}}( 2l +1 )$
can be integrated to $Spin (2l +1)$. $Spin (2l +1)$ and
$\binf$ form a dual pair on $\Flhalf$.
We then obtain a duality theorem for the dual
pair $\Fpairbb$. Indeed it can be stated
literally as Theorem \ref{th_Fpairbb} by simply replacing $\binftwo$
there by our $\binf$ so we will not repeat here.
The dual pair $\Fpairbb$ is related to the dual pair $\Fpairbbtwo$
by choosing a different symmetric bilinear form in defining the
action on $\Flhalf$ of an infinite dimensional Lie subalgebra
of $\hgl$ of $B$ type (cf. Section \ref{subsec_binf}).
\section{Duality in a bosonic Fock space $\Fminusl$}
  \label{sec_Bpair}
\subsection{Untwisted case: dual pair $\Mpaircd$}
     Let us take a pair of bosonic ghost fields
$$ \gamma^+ (z) = \sum_{ n \in {\hz} }
                         \gamma^+_n z^{-n - \hf }, 
\quad \gamma^- (z) = \sum_{ n \in {\hz} }
                         \gamma^-_n z^{-n - \hf} $$
with 
             %
%
the following commutation relations
$$ [ \gamma^+_m, \gamma^-_n ] = \delta_{m+n, 0}, 
\quad [ \gamma^{\pm}_m, \gamma^{\pm}_n ] = 0.$$
We define the Fock space ${\cal F}^{\bigotimes -1}$ 
of the fields $\gamma^+ (z)$ and $\gamma^-(z)$, 
generated by the vacuum $\vac$, satisfying
\begin{eqnarray*}
  \gamma^{+}_n \vac = \gamma^{-}_n \vac = 0
   \quad ( n \in \hn).
\end{eqnarray*}
Now we take $l$ pairs of bosonic ghost fields 
$\gamma^{+,p} (z), \gamma^{-,p} (z) \;(p = 1, \dots, l)$
and consider the
corresponding Fock space $\Fminusl$.

Introduce the following generating functions
\begin{eqnarray}
  E(z,w) & \equiv & \sum_{i,j \in \Z} E_{ij} z^{i -1}w^{-j } 
  = - \sum_{p =1}^l :\gamma^{+,p}  (z) \gamma^{-,p} (w):,
                 \label{eq_hglboson} \\
  e^{pq}_{**} (z) & \equiv &
              \sum_{i,j \in \Z} e^{pq}_{**}(n) z^{ -n -1} 
  = : \gamma^{+,p} (z) \gamma^{+,q} (z): 
    \quad (p \neq q ) ,          \nonumber       \\ 
  e^{pq} (z) & \equiv &
               \sum_{i,j \in \Z} e^{pq}(n)z^{ -n-1}
  = : \gamma^{-,p}(z) \gamma^{-,q} (z):
    \quad (p \neq q ),      \label{eq_genfunc}     \\
  e_*^{pq} (z) & \equiv & \sum_{i,j \in \Z} e_*^{pq} (n)z^{ -n-1}
  = : \gamma^{+,p} (z) \gamma^{-,q} (z):,
   \quad p,q = 1, \dots, l   \nonumber    
\end{eqnarray}
where the normal ordering $::$ means that the operators
annihilating $\vac$ are moved to the right.

It is well known
that the operators $E_{ij}$ $ (i,j \in \Z )$
form a representation in $\Fminusl$ of
the Lie algebra $\hgl$ with central charge $ - l$;
the operators 
 $$ e^{pq}(n), e_*^{pq}(n), 
  e_{**}^{pq}(n) \;( p,q = 1, \dots, l, n \in \Z )$$
form an affine algebra
$\widehat{{\frak {sp}} }(2l)$ with central charge $ -1$ \cite{FF}. 
Denote 
$$ e^{pq} \equiv e^{pq}(0), \; e_*^{pq} \equiv e_*^{pq}(0), 
e_{**}^{pq} \equiv e_{**}^{pq}(0). $$
Then the operators $e^{pq}, e_*^{pq}, e_{**}^{pq}$
$ (p, q = 1, \cdots l )$ form the horizontal subalgebra 
${\frak {sp}}(2l)$ in $\widehat{{\frak {sp}} }(2l)$.
In particular, operators $e_*^{pq}$ $( p, q = 1, \cdots l )$ form
a subalgebra ${\frak {gl}}(l)$ 
in the horizontal subalgebra ${\frak {sp}}(2l)$.
We identify the Borel subalgebra ${\frak b} ({\frak {sp} }(2l) )$
with the one generated by 
$e_{**}^{pq}, e_*^{pq} (p \leq q), \; p, q = 1, \cdots ,l$.

\begin{lemma}
  \begin{enumerate}
  \item[1)]  The action of ${\frak {gl}}(l)$ generated by 
  $ e_*^{pq} \; ( p, q = 1, \cdots, l)$ 
  and that of $\hgl$ generated by 
  $ E_{ij} \; ( i,j \in \Z )$ on $\Fminusl$ commute with each other.
            
  \item[2)]  The action of the horizontal subalgebra ${\frak {sp}}(2l)$
  and that of $\dinf$ generated by 
  $ E_{ij} - E_{1-j, 1-i} \; ( i,j \in \Z )$ on $\Fminusl$ 
 commute with each other.
  \end{enumerate}
\end{lemma}        
Proof is similar to the fermionic case. We omit it here.
 
A similar remark to Remark \ref{rem_grad} holds 
in our bosonic case.
So as a representation of ${\frak {sp}}(2l)$, $\Fminusl$ is decomposed
 into a direct sum of finite dimensional irreducible representations. 
Furthermore as a representation of ${\frak {sp}}(2l)$, 
$\Fminusl$ is isomorphic to the polynomial algebra
$ {\cal P} ( \C^{2l} \bigotimes \C^{\Bbb N} )$,
where ${\frak {sp}}(2l)$ acts on $\C^{2l}\bigotimes \C^{\Bbb N} $ naturally 
on the left of $\C^{2l}$. The action of the Lie algebra
${\frak {sp}}(2l)$ (resp. ${\frak {gl}}(l)$)
can be integrated to an action of $Sp(2l)$ (resp. $GL(l)$). 
The action of $ Sp (2l)$ commutes with the action of $\dinf$ on $\Fminusl$.
Indeed they form a dual pair. Similarly $GL(l)$ and $\hgl$
form a dual pair on $\Fminusl$.
       
Denote by $\Gamma^{+,m}_i $ the $m$-th power of the
determinant of the following $ i \times i$ matrix:
\begin{eqnarray}
\left[ \begin{array}{cccc}
\gamma_{-1/2}^{+,1} & \gamma_{-3/2}^{+,1} & \cdots & \gamma_{-i + 1/2}^{+,1}\\
\gamma_{-1/2}^{+,2} & \gamma_{-3/2}^{+,2} & \cdots & \gamma_{-i + 1/2}^{+,2}\\
\cdots              & \cdots               & \ddots & \cdots                \\
\gamma_{-1/2}^{+,i} & \cdots               & \cdots & \gamma_{-i + 1/2}^{+,i} 
\end{array} \right].
   \label{mat_1}
\end{eqnarray}
We also denote by $\Gamma^{-, m}_{l +1 -i} $ the $m$-th power of
the determinant of the following $ i \times i$ matrix:
\begin{eqnarray}
\left[ \begin{array}{cccc}
\gamma_{- 1/2}^{-,l} & \gamma_{- 3/2}^{-,l} 
& \cdots & \gamma_{-i + 1/2}^{-,l}   \\
\gamma_{- 1/2}^{-,l-1} & \gamma_{- 3/2}^{-,l-1} 
& \cdots & \gamma_{-i + 1/2}^{-,l-1}   \\
\cdots              & \cdots               
& \ddots & \cdots                     \\
\gamma_{- 1/2}^{-,l-i+1} & \cdots               
& \cdots & \gamma_{-i + 1/2}^{-,l-i+1} 
\end{array} \right].
  \label{mat_2}
\end{eqnarray}
We define a map $\Lambda^{\frak{aa} }_{-}: \Sigma (A)  \longrightarrow \hgl_0^*$
by sending $\lambda = (m_1, \cdots, m_l) $ to 
\begin{eqnarray}
  \Lambda^{\frak{aa}}_{-} (\lambda) 
    & = & -m_{i+1} \hL^a_{i-l} + \sum_{k = 1 +i -l}^{-1}
        ( m_{l +k} - m_{l +k +1}) \hL^a_k     \nonumber       \\
    & + & ( m_l -m_1 -l ) \hL^a_0
      + \sum_{k=1}^{i-1} (m_k - m_{k+1} )\hL^a_k + m_i \hL^a_i.
\end{eqnarray}
\begin{theorem}
  \begin{enumerate}
  \item[1)] We have the following $(GL(l), \hgl )$-module
 decomposition:
 \begin{equation}
   \Fminusl = \bigoplus_{\lambda \in \Sigma (A)} I_{\lambda}
       \equiv \bigoplus_{\lambda \in \Sigma (A)} 
          V( GL(l); \lambda) \otimes L \left(\hgl;
                                 \Lambda^{\frak{aa}}_{-}(\lambda), -l
                               \right)
 \end{equation}
 where $V( GL(l); \lambda)$ 
 is the irreducible $GL (l)$-module
 of highest weight $\lambda$ and 
 $L \left(\hgl; \Lambda^{\frak{aa}}_{-}(\lambda), -l \right)$ 
 is the irreducible 
 highest weight module of $\hgl$ with highest weight
 $\Lambda^{\frak{aa}}_{-}(\lambda)$ and central charge $-l$.

 \item[2)] Given $\lambda =(m_1, \ldots, m_l )
 \in \Sigma (A)$, fix $i, j $ such that
 $m_1 \geq \cdots \geq m_i \geq m_{i+1} = \cdots  = m_{j} = 0
  > m_{j+1} \geq \cdots \geq m_l. $
 Then the corresponding highest weight vector in $I_{\lambda}$ 
 with respect to $(GL(l), \hgl )$ is 
  \begin{eqnarray}
   \Gamma^{+, m_1 -m_2}_1 \ldots \Gamma^{+, m_i}_i
    \Gamma^{-,- m_{j+1} }_{j+1} \Gamma^{-, m_{j+1} -m_{j+2} }_{j+2}
   \ldots \Gamma^{-, m_{l-1} - m_l}_l \vac . 
    \label{eq_vec}   
  \end{eqnarray}
 \end{enumerate}
    \label{th_Bpairgl}
 \end{theorem}
\begin{demo}{Proof}
   We first prove that the vector (\ref{eq_vec}) is a highest 
 weight vector for ${\frak {gl}}(l)$. Indeed 
 We see that the action of $e_*^{pq} (p < q, p, q = 1, \ldots, l)$ on 
 the vector(\ref{eq_vec}) has the effect by
 replacing a row in the matrices (\ref{mat_1}) and (\ref{mat_2})
 by another row up in the same matrix, so the determinant of
 the new matrix is zero.
 
 On the other hand, the action of $ E_{i,j} (i \leq 0< j, i, j \in \Z ) $
 on the vector (\ref{eq_vec}) has the effect by
 replacing a column in the matrices (\ref{mat_1}) and (\ref{mat_2})
 by another column with same superscripts but some positive
 subscripts. Then we see that the new entries in the replaced
 column will commute with all other operators in the expression
 of (\ref{eq_vec}), so we can move it over to the right to kill the
 vacuum vector $\vac$. The action of
 $ E_{i,j} (0 <i < j, \mbox{ or } i < j <0, i, j \in \Z ) $
 on the vector(\ref{eq_vec}) has the effect by
 replacing a column in matrices (\ref{mat_1}) and (\ref{mat_2})
 by another column to the left in the same matrix, so 
 the determinant of the new matrix is zero again.
 Thus the vector (\ref{eq_vec}) is also a highest weight vector for $\hgl$.

 One moment's thought shows that the weight 
 of the vector $\Gamma^{+, k}_i \vac $
 with respect to ${\frak h} ({\frak {gl} }(l) )$
 is $ (k, \ldots, k, 0, \ldots, 0)$, where the number of $k$ 
 (resp. $0$) is $i$ (resp. $l -i$).
 It is also easy to see that the weight 
 of the vector $\Gamma^{-, k}_i \vac $
 with respect to ${\frak h} ({\frak {sp} }(2l) )$
 is $ (0, \ldots, 0, -k, \ldots, -k )$, where the number of $-k$ is $i$.
 Then we easily get the highest weight of the vector (\ref{eq_vec})
 is $(m_1, m_2, \ldots, m_l )$.
 
 Easy to see $ E_{a,a} . \Gamma^{-, m}_b = 0 $, for  $a > 0$;
 $ E_{a,a} . \Gamma^{+, m}_b = m \Gamma^{+, m}_b$ if $ 0 \leq a \leq b$,
 $= 0$ if $a > b$.
 Also we have $ E_{a,a} . \Gamma^{+, m}_b = 0 $, for $a \leq 0$;
 $ E_{a,a} . \Gamma^{-, m}_b = - m\Gamma^{-, m}_b$ if $ b-l \leq a \leq 0$,
 $= 0$ if $ -a > b$.
 From these data, we can calculate the 
 following table for the vector (\ref{eq_vec}):

 \vspace{0.1in}
 \begin{center}
 \begin{tabular}{|l|l|l|l|l||l|l|l|l|} \hline
  $k$      & $i$   & $i-1$     & \ldots & 1     & 0   
   & $-1 $ & \ldots & $j-l+1$     \\ \hline
 $ E_{kk}$ & $m_i$ & $m_{i-1}$ & \ldots & $m_1$ & $m_l$ 
   & $m_{l-1}$ & \ldots &  $m_{j+1}$     \\ \hline
 \end{tabular}
 \end{center}
 \vspace{0.1in}

 Now we can easily see that the highest weight
 of the vector (\ref{eq_vec}) is indeed given by 
 $ \Lambda^{\frak{aa} }_{-}(\lambda)$.

  Since $\lambda $ ranges over all $\Sigma (A)$, the 
 decomposition in 1) follows by general nonsense of dual pair theory.
\end{demo}

\begin{remark}
  The first part of the theorem was proved in \cite{KR2}. The most difficult
 part of their proof is to determine the highest weight 
 $\hL^{\frak{aa} }_{} (\lambda)$ 
 for the vector (\ref{eq_vec}) is rather indirect and based on
 a series of combinatorial lemmas.
 As we see from our proof, this follows from the explicit formula
 of the highest weight vector (\ref{eq_vec}) fairly easily.
\end{remark}

We now define a map $\Lambda^{\frak{cd} }$ from $\Sigma (C)$ to ${\dinf}_0^*$ 
which maps $\lambda = (m_1, \ldots, m_l) $ to
$$ \Lambda^{\frak{cd} } (\lambda) =
    (-2l -m_1 -m_2) \hL^d_0 + \sum_{k =1}^l (m_k -m_{k +1}) \hL^d_k $$
with the convention here and below $ m_{l+1} = 0.$
\begin{theorem}
 \begin{enumerate}
   \item[1)] We have the following $\Mpaircd$-module
 decomposition:
 \begin{equation}
   \Fminusl = \bigoplus_{\lambda \in \Sigma (C)} I_{\lambda}
       \equiv \bigoplus_{\lambda \in \Sigma (C)} 
          V( Sp(2l); \lambda) \otimes L \left(\dinf;
                                 \Lambda^{ \frak{cd}} (\lambda), -l
                               \right)
 \end{equation}
 where $V( Sp(2l); \lambda)$ is the irreducible $Sp (2l)$-module
 of highest weight $\lambda$ and 
 $L \left(\dinf; \Lambda^{ \frak{cd}} (\lambda), -l  \right)$ 
 the irreducible highest weight $\dinf$-module of highest weight
 $\Lambda^{ \frak{cd}} (\lambda) $
 and  central charge $-l$.

 \item[2)] Given $\lambda = (m_1, \ldots, m_l) \in \Sigma (C)$,
 with respect to $\Mpaircd$
 the isotypic subspace $I_{\lambda}$ is irreducible with
 highest weight vector
 \begin{eqnarray}
  \Gamma^{+, m_1- m_2}_1 \Gamma^{+, m_2- m_3}_2 
    \ldots \Gamma^{+, m_l-m_{l+1}}_l \vac .     
    \label{vec_bos}
 \end{eqnarray}
 \end{enumerate}
  \label{th_Bpair}
 \end{theorem}
\begin{demo}{Proof}
  Since the vector (\ref{vec_bos}) is a highest weight vector
 for $\hgl$, it is so for the Lie subalgebra $\dinf$ of $\hgl$.
 On the other hand, $e_{**}^{pq} (p, q = 1, \ldots, l)$ 
 annihilates the vector (\ref{vec_bos}) since it commutes with
 all $\Gamma^{+, k}_i \; (i = 1, \ldots, l ) $
 and it annihilates $\vac$. Thus the vector (\ref{vec_bos}) 
 is also a highest weight vector for $Sp(2l)$.

 The highest weight of the vector (\ref{vec_bos}) for $Sp(2l)$ is 
 already calculated as in Theorem \ref{th_Bpairgl} since
 ${\frak {sp}}( 2l)$ and  ${\frak {gl}}( l)$ share the same
 Cartan subalgebra. The highest weight of the vector (\ref{vec_bos}) 
 with respect to $\dinf$ can be seen from the following table:

 \vspace{0.1in}
 \begin{center}
 \begin{tabular}{|l|l|l|l|l|} \hline
  $k$      & $i$   & $i-1$     & \ldots & 1
    \\ \hline
 $ E_{kk}$ & $m_i$ & $m_{i-1}$ & \ldots & $m_1$ 
    \\ \hline
 \end{tabular}
 \end{center}
 \vspace{0.1in}

 Then we can read off the highest weight of
 the vector (\ref{vec_bos}) to be $\Lambda^{ \frak{cd}}(\lambda)$.

 The decomposition in 1) now follows from general nonsense
 of dual pair theory since all irreducible representations 
 of $Sp(l)$ already appear in the decomposition.
\end{demo}
We immediately have the following corollary.
\begin{corollary}
    The space of invariants of $Sp (2l)$ in the Fock space
  $\Fminusl$ is the irreducible
  module $L(\dinf; -2l \hL^d_0 )$ of central charge $ -l$.
\end{corollary}
\subsection{Twisted case: dual pairs $\Mpairdc$}
It is known \cite{FF} that
the Fourier components of the following generating functions
\begin{eqnarray}
  e^{pq}_{**} (z) & \equiv & \sum_{i,j \in \Z} e^{pq}_{**}(n) 
                   z^{ -n -1}
  = : \gamma^{+,p} (z) \gamma^{+,q} ( -z): 
    \quad (p \neq q )           \nonumber       \\ 
  e^{pq} (z) & \equiv & \sum_{i,j \in \Z} e^{pq}(n) 
                   z^{ -n -1}
  = : \gamma^{-,p}(z) \gamma^{-,q} ( -z):
    \quad (p \neq q )      \label{eq_series}       \\
  e_*^{pq} (z) & \equiv & \sum_{i,j \in \Z} e_*^{pq} (n) 
                   z^{ -n -1}
  = : \gamma^{+,p} (z) \gamma^{-,q} (z): 
   \quad (p,q = 1, \dots, l  ) \nonumber    
\end{eqnarray}
span an affine algebra ${\frak {gl}}^{(2)} (2l)$
of type $A^{(2)}_{2l -1}$ of central charge $ -1$ when acting on $\Fminusl$. 
The horizontal subalgebra of the affine algebra ${\frak {gl}}^{(2)} (2l)$
spanned by $e^{pq}_{**} (0), e_*^{pq} (0),$ $ e^{pq} (0) $
$(p, q = 1, \ldots, l)$ 
is isomorphic to the Lie algebra ${ \frak {so}} (2l)$.
Put 
\begin{eqnarray*}
 \sum_{i,j \in \Z} (E_{ij} - ( -1)^{i+j} E_{1-j, 1-i})z^{i-1}w^{-j}
  = \sum_{ k=1}^l \left(:\gamma^{+,k} (z) \gamma^{-,k} (w):
                        + :\gamma^{+,k} ( -w) \gamma^{-,k} ( -z):
                  \right).
\end{eqnarray*}
Proof of the following lemma is straightforward.
\begin{lemma}
    The action of the horizontal subalgebra ${\frak {so}}(2l)$
  and that of $\cinf $ generated by the operators
  $ E_{ij} - ( -1)^{i+j} E_{1-j, 1-i} \; ( i,j \in \Z )$ on $\Fminusl$ 
 commute with each other. 
\end{lemma}        

The action of ${\frak {so}}(2l)$ can be integrated to an action
of the Lie group $SO (2l)$ and extends to $O (2l)$ naturally.
Indeed $O(2l)$ and $\cinf$ form a dual pair on $\Fminusl$.

Define $\Gamma_j^{det} $ to be the determinant of the following
$ (2l -j) \times (2l -j)$ matrix $M$:
\begin{eqnarray*}
 \left[ \begin{array}{ccccccc}
\gamma_{-1/2}^{+,1}& \gamma_{-1/2}^{+,2}& \cdots& \gamma_{-1/2}^{+,l} &
 \gamma_{- 1/2}^{-,j+1} &  \cdots     &    \gamma_{- 1/2}^{-,l}    \\
\gamma_{-3/2}^{+,1}& \gamma_{-3/2}^{+,2}& \cdots& \cdots  &
-\gamma_{-3/2}^{-,j+1}  &  \cdots     &    -\gamma_{-3/2}^{-,l}       \\
\cdots             & \cdots             & \ddots& \cdots              &
\cdots                 & \cdots      & \cdots             \\
\gamma_{j-2l+1/2}^{+,1}& \gamma_{j-2l + 1/2}^{+,2} & \cdots& \gamma_{j-2l+1/2}^{+,l} & 
(-1)^{j+1}\gamma_{j-2l+1/2}^{-,j+1} & \cdots &
(-1)^{j+1}\gamma_{j-2l+1/2}^{-,l}
\end{array} \right].
\end{eqnarray*}
Note that in the last $l-j$ columns there is a sign $ (-1)^{i+1}$
in the $i$-th row.

Denote by $\widetilde{\Gamma}_l^{+, m }$ the $m$-th power of the
deteminant of the matrix obtained from the matrix (\ref{mat_1}) 
for $i =l$ by replacing the last row with the vector 
$$ \left( \gamma_{-1/2}^{- ,i}, -\gamma_{-3/2}^{- ,i}, \ldots,
    (-1)^{i+1} \gamma_{-i + 1/2}^{-,i}
   \right). $$

Define a map $\Lambda^{ \frak{dc}}$
from $\Sigma(D)$ to ${\cinf}_)^*$
by sending 
$ \lambda = (m_1, \cdots, \overline{ m}_l) \; (m_l > 0 )$ to 
$$
 \Lambda^{ \frak{dc}} (\lambda) =
    (-l -m_1) \hL_0^c + \sum_{k =1}^l (m_k -m_{k+1}) \hL_k^c, 
$$ 
sending $(m_1, \cdots, m_j, 0, \cdots, 0 ) \;(j<l)  $ to 
$$
 \Lambda^{ \frak{dc}} (\lambda) =
    (-l -m_1) \hL_0^c + \sum_{k =1}^j (m_k -m_{k+1}) \hL_k^c, 
$$
and sending $(m_1, \cdots, m_j, 0, \cdots, 0 ) \bigotimes {det} \;(j<l)$ to 
$$ 
 \Lambda^{ \frak{dc}} (\lambda) =
   (-l -m_1) \hL_0^c + \sum_{k =1}^{j-1} (m_k -m_{k+1}) \hL_k^c
     + (m_j -1) \hL_j^c + \hL_{2l-j}^c
$$
if $m_1 \geq \ldots  m_{j} > m_{j+1} = \ldots = m_l = 0.$
\begin{theorem}
 \begin{enumerate}
   \item[1)] We have the following $(O(2l), \cinf)$-module
 decomposition:
 \begin{equation}
   \Fminusl = \bigoplus_{\lambda \in \Sigma (D)} I_{\lambda}
       \equiv \bigoplus_{\lambda \in \Sigma (D)} 
          V( O(2l); \lambda) \otimes L \left(\cinf;
                                 \Lambda^{ \frak{dc}} (\lambda), -l
                               \right)
 \end{equation}
 where $V( O(2l); \lambda)$ is the irreducible $O (2l)$-module
 parametrized by $\lambda \in \Sigma (D)$ and 
 $L \left(\cinf; \Lambda^{ \frak{dc}} (\lambda) , -l \right)$
 is the irreducible $\cinf$-module of highest weight
 $\Lambda^{ \frak{dc}} (\lambda) $
 and central charge $ -l$.
  \item[2)] With respect to $ ( {\frak {so} } (2l), \cinf )$, 
  \begin{enumerate}
  \item[a)] 
   the isotypic subspace $ I_{\lambda} $
  is decomposed into a sum of two irreducible representations
  with highest weight vectors  
  \begin{equation}
  \Gamma_1^{+, m_1 -m_2} \cdots \Gamma_{l-1}^{+, m_{l-1} -m_l }
               \Gamma_l^{+, m_l }  \vac  
    \label{eq_hwt1}
 \end{equation}
  and
 \begin{equation}
  \Gamma_1^{+, m_1 -m_2} \cdots 
      \Gamma_{l-1}^{+, m_{l-1} - m_l} \widetilde{\Gamma}_l^{+, m_l } \vac 
    \label{eq_hwt2}
 \end{equation}
 in the case $\lambda = (m_1, \ldots, \overline{ m}_l )
 \in \Sigma(D), m_l > 0.$ 
 The highest weight of (\ref{eq_hwt1}) for ${\frak {so} } (2l)$
 is $(m_1, \ldots, m_{l-1}, m_l)$
 while that of (\ref{eq_hwt2}) for ${\frak {so} } (2l)$
 is $(m_1, \ldots, m_{l-1}, - m_l)$;
 \item[b)] 
  the isotypic subspace $ I_{\lambda} $ is irreducible
  with highest weight vector
  \begin{equation}
   \Gamma_1^{+, m_1 -m_2 } \cdots \Gamma_{l-1}^{+, m_{l-1} -m_l }  \vac  
        \label{eq_hwt3}
 \end{equation}
 in the case $ \lambda = (m_1, \cdots, m_{l-1}, 0 ) \in \Sigma(D)$
 with $m_{l} = 0; $
  \item[c)] 
  the isotypic subspace $ I_{\lambda} $ is irreducible
  with highest weight vector
  \begin{eqnarray}
   \Gamma_1^{+, m_1-m_2 } \cdots \Gamma_j^{+, m_{j-1} -m_j }
     \Gamma_j^{+, m_j -1 } \Gamma_j^{det} \vac 
  \label{eq_hwt4}
 \end{eqnarray}
 in the case $ \lambda = (m_1, \cdots, m_j, 0, \ldots, 0 ) \bigotimes {det}
  \in \Sigma(D) ,$ 
 $ m_1 \geq \dots \geq m_j > 0, m_{j+1} = \ldots = m_l = 0, 0 \leq j < l. $
   \end{enumerate}
 \end{enumerate}
   \label{th_Mpairdc}
\end{theorem}
\begin{demo}{Proof}
   One can prove part 2a) and 2b) in a similar way as in the proofs of
 Theorems \ref{th_Bpairgl} and \ref{th_Bpair}.  In the case of 2c), 
 it suffices to check that $\Gamma_0^{det} \vac $ is a highest weight
 vector for ${\frak {so} } (2l)$ since the highest weight vectors in
 Bosonic Fock space form a semi-group. Indeed an action of the
 element $\sum_{ m \in \hz} :\gamma_{-m}^{+,p} \gamma_{m}^{-,q} \; ( p <q)$
 (resp. $\sum_{ m \in \hz} :\gamma_{-m}^{+,p} \gamma_{m}^{+,q}$)
 in the Borel subalgebra of ${\frak {so} } (2l)$ on 
 $\Gamma_0^{det} \vac $ has the effect of replacing the $q$-th 
 (resp. ($l+q$)-th)
 column of the matrix $M$ by
 its $p$-th column. So the determinant of the matrix thus obtained
 is zero.

 Consider the element 
 $$ E_{i,j} - ( -1)^{i+j} E_{1-j, 1-i} 
   \sum_{p =1}^l :\gamma_{-i+1/2}^{+,p} \gamma_{j - 1/2}^{-,p}: 
    - ( -1)^{i+j}:\gamma_{j - 1/2}^{+,p}  \gamma_{-i+1/2}^{-,p}:
 $$
 in the Borel of $\cinf$, where $i\leq 0< j, i, j \in \Z $.
 By applying $ :\gamma_{-i+1/2}^{+,p} \gamma_{j - 1/2}^{-,p}: $
 to $\Gamma_0^{det} \vac $, we obtain up to a sign the
 determinant of the  $ (2l -j-2) \times (2l -j-2)$ submatrix of $M$ by
 deleting the $(-i +1) $-th and $j$-th rows and $p$-th and $(l+p) $-th
 column. The sign can be determined to be 
 $ - ( -1)^i (-1)^{j + p + (-i+1) + (l+p-1)} = (-1)^{j +l +1} . $
 By applying $:\gamma_{j - 1/2}^{+,p}  \gamma_{-i+1/2}^{-,p}:$
 to $\Gamma_0^{det} \vac $, we obtain up to a sign the
 determinant of the same submatrix of $M$. The sign can be determined
 to be $ - ( -1)^{j-1} (-1)^{(-i+1) + p + (j-1)+ (l+p-1)} = (-1)^{i +l +1} . $
 Thus $ E_{i,j} - ( -1)^{i+j} E_{1-j, 1-i} $ acting on
 $\Gamma_0^{det} \vac $ is zero due to cancellations.

 On the other hand,
 $ E_{i,j} - ( -1)^{i+j} E_{1-j, 1-i}\; (0< i < j)  $
 acts on $\Gamma_0^{det} \vac $ has the effect by 
 replacing the $j$-th row by $i$-th row in $M$, so the determinant
 of the matrix thus obtained is obvious zero. The case for $i < j <0$
 is similar.

 One can determine the highest weights for the vectors
 in 2a) and 2b) with respect to ${\frak {so} } (2l)$ and $\cinf$
 in a similar way as in the proofs of
 Theorems \ref{th_Bpairgl} and \ref{th_Bpair}.  
 A simple calculation shows that
 the highest weights of vectors in 2) with respect to 
 ${\frak {so} } (2l)$ are indeed given by $\lambda \in \Sigma(D)$.
 The highest weight of $\Gamma_j^{det} \vac$ with respect to 
 ${\frak {so} } (2l)$ is $(1, \ldots, 1, 0, \ldots, 0)$ where
 the number of $1$ is equal to $j$.
 The highest weight of the vector (\ref{eq_hwt1}) and (\ref{eq_hwt3}) 
 with respect to $\cinf$ can be read off from the following table
 which is equal to $\Lambda^{ \frak{dc}} (\lambda)$: 
 
 \vspace{0.1in}
 \begin{center}
 \begin{tabular}{|l|l|l|l|l||l|l|l|l|} \hline
  $k$      & $i$   & $i-1$     & \ldots & 1     & 0   
   & $-1 $ & \ldots & $j-l+1$     \\ \hline
 $ E_{kk}$ & $m_i$ & $m_{i-1}$ & \ldots & $m_1$ & $m_l$ 
   & $m_{l-1}$ & \ldots &  $m_{j+1}$     \\ \hline
 \end{tabular}
 \end{center}
 \vspace{0.1in}

 The vector (\ref{eq_hwt2}) is obtained from (\ref{eq_hwt1}) 
 by an action of the element $\tau$ of $O(2l)$ as defined
 in \ref{eq_tau}. So it has the same highest weight with respect to 
 $\cinf$ since $O(2l)$ commutes with $\cinf$. The highest
 weight of $\Gamma_j^{det} \vac$ can be read off 
 to be $(-l -1) \hL^c_0 + \hL_{2l-j}$ from the following table: 

 \vspace{0.1in}
 \begin{center}
 \begin{tabular}{|l|l|lll|} \hline
  $k$                 & $<2l-j$ & $2l-j,$ & \ldots, & 1
    \\ \hline
 $E_{kk}-E_{1-k,1-k}$ &   0     & 1,      & \ldots, & 1 
    \\ \hline
 \end{tabular}
 \end{center}
 \vspace{0.1in}

 The highest weights with respect to $\cinf$ form a semigroup.
 From this we can see that the highest weight of (\ref{eq_hwt4}) is
 indeed given as in the theorem.
\end{demo}
\section{Duality in the Fock space $\Fminuslpmhalf$ }
  \label{sec_super}
\subsection{Dual pair $\Mpairospd$}
 We denote by $\Fhalfminusl$ the tensor product
of the Fock space $\Fminusl$ of $l$ pairs of bosonic ghost fields
and the Fock space $\Fhalf$ of a neutral fermionic field.
It is known \cite{FF} that the Fourier components of 
the generating functions $e^{pq}(z), e^{pq}_* (z),$ 
$ e^{pq}_{**}(z)$
in (\ref{eq_genfunc}) and the following generating functions
\begin{eqnarray}
 e^{p} (z) \equiv \sum_{n \in \Z}
 e^{p}(n) \NN & = & : \gamma^{-,p}(z) \phi (z):, 
         \nonumber       \\
 \tilde{e}_*^{p} (z) \equiv \sum_{n \in \Z}
 \tilde{e}_*^{p}(n) \NN & = & : \gamma^{+,p}(z) \phi (z):, 
    \quad p = 1, \dots, l      \nonumber 
  \label{eq_oddsuper}
\end{eqnarray}
span the affine superalgebra 
$\widehat{ {\frak {osp}} }(1, 2l)$ of level $-1$. Put
\begin{eqnarray*}
e^{pq} \equiv e^{pq}(0), & e_*^{pq} \equiv e_*^{pq}(0), &
e^{p} \equiv e^{p} (0),          \\
e_{**}^{pq} \equiv e_{**}^{pq}(0),&
 e_*^{p} \equiv e_*^{p}(0), & p, q = 1, \dots, l
\end{eqnarray*}
Then the operators $ e^{pq}, e_*^{pq}, e_{**}^{pq},
\tilde{e}^{p}, \tilde{e}_*^{p}$ $ (p, q = 1, \dots, l) $ 
generate the horizontal subalgebra
${\frak {osp}}(1,2l)$ of the affine superalgebra
$\widehat{{\frak {osp}} }(1, 2l) $.
We identify the Borel subalgebra ${\frak b} ({\frak {osp} }(1, 2l) )$
with the one generated by 
$  e_*^{pq} (p \leq q), e_{**}^{pq}, \tilde{e}_*^{p},$ 
$ p, q = 1, \cdots, l. $

The following lemma can be proved similarly 
as in the fermionic Fock space case, cf. Lemma \ref{lem_oddabel}.
\begin{lemma}
  Putting
 \begin{equation}
      \begin{array}{rcl}
      \lefteqn{
  \sum_{i,j \in \Z}  \left( E_{ij}  - E_{1-j,1-i}
                     \right) z^{i- 1} w^{-j}
              }                                \\
       & &
  = \sum^l_{k=1} \left(: \gamma^{+,k} (z) \gamma^{-,k} (w) :
                     - : \gamma^{+,k} (w) \gamma^{-,k} (z) :
                 \right) 
                     + :\phi (z) \phi (w):
      \end{array}
 \end{equation}
  defines a representation of $\dinf$ of central
  charge $-l + \hf$ on $\Fhalfminusl$.
\end{lemma}
The following lemma can also be proved similarly as in the fermionic 
Fock space case.
\begin{lemma}
  The action of the Lie supergroup $Osp(1, 2l)$ by integrating the 
 horizontal subalgebra ${\frak {osp}}(1, 2l)$ commute with 
 that of Lie algebra $\dinf$ generated by 
  $ E_{i,j} - E_{1-j, 1-i} ( i,j \in \Z )$ on $\Fhalfminusl$.
\end{lemma}        

$Osp(1, 2l)$ and $\dinf$ form a dual pair on $\Fhalfminusl$.
We define a map $\Lambda^{ \frak{ospd}}$
from $\Sigma (Osp) $ to ${\dinf}^*_0$ by sending
$ \lambda = (m_1, \ldots, m_l )$ to
$$ \Lambda^{ \frak{ospd}} (\lambda )=
   (-2l +1 -m_1 -m_2) \hL^d_0 + \sum_{k =1}^l (m_k - m_{k +1} )\hL^d_k.$$
We obtain the following duality theorem.
\begin{theorem}
 \begin{enumerate}
  \item[1)]  We have the following $ \Mpairospd$-module decomposition
   \begin{eqnarray*}
    \Fhalfminusl = \bigoplus_{\lambda \in \Sigma (Osp)} 
          V( Osp(1, 2l); \lambda) \otimes L \left(\dinf;
                             \Lambda^{ \frak{ospd}} (\lambda), -l +1/2 
                               \right)
   \end{eqnarray*}
 where $V( Osp; \lambda)$ is the irreducible
 module of $Osp (1, 2l)$ of highest weight $\lambda$, and 
 $L \left(\dinf; \Lambda^{ \frak{ospd}} (\lambda), -l +1/2 \right)$
 is the irreducible 
 highest weight $\dinf$-module of highest weight
 $\Lambda^ \frak{ospd} (\lambda) $
 and central charge $-l +1/2$.

 \item[2)] Given $\lambda = (m_1, \ldots, m_l) \in \Sigma (Osp)$,
 the isotypic subspace $I_{\lambda}$ is irreducible
 with highest weight vector
 $$ \Gamma^{+, m_1- m_2}_1 \ldots \Gamma^{+, m_l-m_{l+1}}_l \vac .     $$
 \end{enumerate}
\end{theorem}
\subsection{Twisted case: dual pair $\Mpairbc$}
   Recall that a bosonic field
 $ \chi (z) = \sum_{ n \in \hz} \chi_n z^{ -n - \hf}$ 
satisfies the following commutation relations:
$$[ \chi_m , \chi_n ] = ( -1)^{m +\hf} \delta_{m, -n}, \quad m,n \in \hz .$$
Let $\Fminuslhalf$ be the tensor product of the Fock space 
of $l$ pairs of bosonic ghost fields
$\gamma^{\pm, k} (z) \; ( k = 1, \ldots, l)$
and the Fock space $\Fminushalf$ of $\chi (z)$.

It is known \cite{FF} that the Fourier components of the generating functions 
$ \tilde{e}^{pq} (z), e^{pq}_* (z), \tilde{e}^{pq}_{**} (z)$ 
as in (\ref{eq_series}) and the following generating functions
\begin{eqnarray*}
 \zeta (z) \equiv \sum_{n \in \Z}
            \zeta(n) \NN & = & : \chi (z) \chi ( -z):,   \nonumber     \\
 \tilde{e}^{p} (z) \equiv \sum_{n \in \Z}
            e^{p}(n) \NN & = & :
              \gamma^{-,p}(z) \chi ( -z):,   \nonumber     \\
 \tilde{e}_*^{p} (z) \equiv \sum_{n \in \Z}
           e_*^{p}(n) \NN & = & : \gamma^{+,p}(z) \chi (z):  \nonumber    \\
\end{eqnarray*}
span an affine algebra $A^{(2)}_{2l}$ of central charge $-1$. 
Its horizontal subalgebra is isomorphic to ${\frak {so}}(2l +1)$.
Let 
\begin{eqnarray*}
 \begin{array}{rcl}
  \lefteqn{
   \sum_{i,j \in \Z} (E_{ij} - ( -1)^{i +j} E_{1-j, 1-i})z^{ i-1} w^{-j}
          }              \\
  & = & \sum_{ k=1}^l \left( :\gamma^{+,k} (z) \gamma^{-,k} (w):
                             + :\gamma^{+,k} ( -w) \gamma^{-,k} ( -z):
                      \right)
              + : \chi (z) \chi (-w):.
  \end{array}
\end{eqnarray*}
\begin{lemma}
    The action of the horizontal subalgebra ${\frak {so}}(2l +1)$
  and that of $\cinf $ generated by 
  $ E_{ij} - ( -1)^{i +j}E_{1-j, 1-i} \; ( i,j \in \Z )$ on $\Fminuslhalf$ 
 commute with each other.
\end{lemma}        
The action of ${\frak {so}}(2l +1)$ can be lifted to $O(2l +1)$.
$ O(2l +1)$ and $\cinf$ form a dual pair on $\Fminuslhalf$.
Define $\widetilde{\Gamma}_j^{det} $ to be the determinant of the following
$ (2l -j +1) \times (2l -j +1)$ matrix $\widetilde{M}$:
\begin{eqnarray*}
 \left[ \begin{array}{ccccccc}
\gamma_{-1/2}^{+,1}& \cdots & \gamma_{-1/2}^{+,l} &
 \gamma_{- 1/2}^{-,j+1} &  \cdots & \gamma_{- 1/2}^{-,l} & \chi_{- 1/2}   \\
\gamma_{-3/2}^{+,1}& \cdots & \cdots  &
-\gamma_{-3/2}^{-,j+1}  &  \cdots & -\gamma_{-3/2}^{-,l} & -\chi_{- 3/2}   \\
\cdots             & \cdots & \cdots              &
\cdots                 & \cdots  & \cdots  & \cdots             \\
\gamma_{j-2l-1/2}^{+,1}& \cdots& \gamma_{j-2l-1/2}^{+,l} & 
(-1)^j\gamma_{j-2l-1/2}^{-,j+1} & \cdots &
(-1)^j \gamma_{j-2l-1/2}^{-,l} & (-1)^j \chi_{-2l +j - 1/2}
\end{array} \right].
\end{eqnarray*}
Define a map $ \Lambda^{ \frak{bc}}$ from $\Sigma(B)$ to ${\cinf}_0^*$ 
by sending
$$\lambda = (m_1, m_2, \ldots, m_l)$$
 to
$$ \Lambda^{ \frak{bc}} (\lambda) =
     (-l -m_1 -1/2) \hL_0^c + \sum_{k =1}^j (m_k -m_{k+1}) \hL_k^c $$
and sending
$$\lambda = (m_1, m_2, \ldots, m_l) \bigotimes det $$
 to
$$ \Lambda^{ \frak{bc}} (\lambda) =
     (-l -m_1 -1/2) \hL_0^c + \sum_{k =1}^{j-1} (m_k -m_{k+1}) \hL_k^c
    +(m_j -1) \hL^c_j + \hL^c_{2l -j +1} $$
where $ m_1 \geq \ldots \geq m_j > m_{j+1} = \ldots = m_l = 0 .$
Proof of the following theorem is similar to that of Theorem~\ref{th_Mpairdc}.
\begin{theorem}
 \begin{enumerate}
   \item[1)] We have the following $( O(2l +1), \cinf)$-module
 decomposition:
 \begin{eqnarray*}
   \Fminuslhalf = \bigoplus_{\lambda \in \Sigma (B)} I_{\lambda}
       \equiv \bigoplus_{\lambda \in \Sigma (B)} 
          V( O(2l +1); \lambda) \otimes L \left(\cinf;
                                 \Lambda^{ \frak{bc}} (\lambda), -l -1/2
                               \right)
 \end{eqnarray*}
 where $V( O(2l +1); \lambda)$ is the irreducible $O(2l +1)$-module
 parametrized by $\lambda \in \Sigma (B)$ and 
 $L \left(\cinf; \Lambda^{ \frak{bc}} (\lambda) -l -1/2 \right)$ 
is the irreducible highest weight $\cinf$-module of highest weight
 $\Lambda^{ \frak{bc}} (\lambda) $
 and central charge $-l -1/2$.

 \item[2)] With respect to ${\frak{so}} (2l +1)$, 
 the isotypic subspace $I_{\lambda}$ is decomposed into
 a sum of two irreducible representations with highest weight vectors
  \begin{eqnarray*}
   \Gamma^{+, m_1- m_2}_1 \ldots
           \Gamma^{+, m_{j-1} -m_j}_{j-1} \Gamma^{+, m_j}_j \vac
  \mbox{ for } \lambda =(m_1, m_2, \ldots, m_l)  
  \end{eqnarray*}
  and 
  \begin{eqnarray*}
   \Gamma^{+, m_1- m_2}_1 \ldots
        \Gamma^{+, m_{j-1} -m_j}_{j-1} \Gamma^{+, m_j -1}_j 
            \widetilde{\Gamma}_j^{det} \vac
  \mbox{ for } \lambda =(m_1, m_2, \ldots, m_l) \bigotimes det
  \end{eqnarray*}
  if $ m_1 \geq \ldots \geq m_j > m_{j+1} = \ldots = m_l = 0 .$
 \end{enumerate}
\end{theorem}
\section{Reciprocity laws and tensor categories of
modules over $\binf$, $\binftwo$, $\cinf$ and $\dinf$}
  \label{sec_reciprocity}
In this section we outline an application of the theory of
dual pairs we have developed in this paper, namely we establish 
certain reciprocity laws as a formal concequence of the so-called
``see-saw'' pairs (cf. \cite{H2} and references therein).
There will be a number of reciprocity laws associated to
different (series of) dual pairs given in the two tables
in the introduction.
We will take a step further to establish an equivalence
of two tensor categories as in \cite{W}. One tensor category consists of
modules over finite dimensional Lie groups of a fixed type with
all ranks. The tensor product
in this category is defined in terms of an induction functor.
The other consists of certain modules over
an infinite dimensional Lie algebra with the usual tensor product. 
This formulation of equivalence of tensor categories is
particularly natural in our infinite dimensional setting.
\subsection{Reciprocity laws associated to see-saw pairs}
   We formulate in this subsection a theory of
reciprocity laws associated to see-saw pairs in a general term
following the treatment of \cite{H2}.

Given two dual pairs $(G, {\frak g}' )$ and $ (K, {\frak k}')$ 
acting on a same vector space $V$. For our purpose we
assume that $G$ and $K$ are compact (or complex) Lie groups
and the irreducible representations of $G$ and $K$ are all
finite dimensional. Assume that these two dual pairs
satisfy the following relations: $ K \subset G$ and
$ {\frak g}' \subset {\frak k}'$. Such a pair of dual pairs
is called a see-saw pair. 

A representation $\rho$ of $G$ can be decomposed into a direct 
sum of irreducible $H$-modules $\sigma$ with multiplicities
$m_{ \rho \sigma}$. Denote
$\rho \mid_K = \sum_{\sigma} m_{ \rho \sigma} \sigma.$
On the other hand a representation of ${\frak k}'$ decomposes
into a sum of irreducible ${\frak g}'$-modules:
$ \sigma' \mid_{{\frak g}'} \cong
   m^{'}_{ {\sigma}' {\rho}'} {\rho}'.$

Now we decompose $V$ as a $K \times {\frak g}'$-module.
Let $\sigma$ be an irreducible module of $K$ and
${\rho}'$ an irreducible module of $ {\frak g}'$. We compute
the $\sigma \bigotimes {\rho}'$-isotypic component
$V^{\sigma \bigotimes {\rho}'}$ in two different ways
and then compare the results. We have
$$ V^{\sigma \bigotimes {\rho}'} = (V^{\sigma})^{{\rho}'}
 \cong (\sigma \bigotimes {\sigma}')^{{\rho}'}
 \cong \sigma \bigotimes ({{\sigma}'}^{{\rho}'})
 \cong \mbox{Hom}_{{\frak g}'} ({\rho}', {\sigma}') 
            \bigotimes ({\sigma} \bigotimes {\rho}')
$$
and
$$ V^{\sigma \bigotimes {\rho}'} = (V^{{\rho}'})^{\sigma}
 \cong ( \rho \bigotimes {\rho}')^{\sigma}
 \cong \rho^{\sigma} \bigotimes {\rho}'
 \cong \mbox{Hom}_{K} (\sigma, \rho) \bigotimes ({\sigma} \bigotimes {\rho}').
$$
It follows that
\begin{eqnarray*}
  \mbox{Hom}_{K} (\sigma, \rho)
  \cong \mbox{Hom}_{K \times {\frak g}'}({\sigma} \bigotimes {\rho}', V)
  \cong \mbox{Hom}_{{\frak g}'} ({\rho}', {\sigma}').
\end{eqnarray*}
In particular we have 
\begin{equation}
 m^{'}_{ {\sigma}' {\rho}'} {\rho}' = m_{ \rho \sigma}.
 \label{rec_general}
\end{equation}

\subsection{Reciprocity laws and equivalence of tensor categories}
 We use the dual pair $\Fpair$ acting on $\Fl$ to demonstrate 
how the general reciprocity law applies to our situation. To emphasize the
rank $l$, we will write $\Sigma(D_l)$ for $\Sigma(D)$. 
Then by Theorem \ref{th_Fpair} we have the following decompositions:
 \begin{eqnarray}
   \Fm & = &  \bigoplus_{\mu \in \Sigma (D) } 
          V(O(2m); \mu) \otimes L \left(\dinf;
                                 \Lambda^{\frak{dd}} (\mu), m
                               \right)
   \label{decom_Fpairm}    \\
   \Fn & = &  \bigoplus_{\nu \in \Sigma (D) } 
          V(O(2n); \nu) \otimes L \left(\dinf;
                                 \Lambda^{\frak{dd}} (\nu), n
                               \right)
   \label{decom_Fpairn}    \\
   \Fmn & = &  \bigoplus_{\lambda \in \Sigma (D) } 
          V(O(2m +2n); \lambda) \otimes L \left(\dinf;
                                 \Lambda^{\frak{dd}} (\lambda), m+n
                               \right).
   \label{decom_Fpairmn}
 \end{eqnarray}
By the obvious isomophism 
$ \Fmn \cong \Fm \bigotimes \Fn$, we see that 
$\Fpairmtimesn$ also form a dual pair on $\Fmn$. The two dual pairs 
\begin{eqnarray*}
 \left\{ \begin{array}{cc}
  \Bigl( O(2m +2n), 
    & 
    \dinf |_{c=m +n } \Bigl) \nonumber \\
 \uparrow                          & \downarrow \nonumber\\
 \Bigl( O(2m)\times O(2n), 
     & \dinf |_{c=m } \bigoplus \dinf |_{c=n } \Bigl)
                       \end{array} \right. 
\nonumber
\end{eqnarray*}
form a see-saw pair, where inclusions of Lie groups/algebras 
are shown by the arrows and the second inclusion is given
by the diagonal imbedding. So as a consequence of
(\ref{rec_general}), we have
\begin{eqnarray}
 m \left( \Lambda^{\frak{dd}} (\mu) \otimes \Lambda^{\frak{dd}} (\nu),
        \Lambda^{\frak{dd}} (\lambda)
      \right)
 = m ( \lambda\mid_{ O(m) \times O(n)}, \mu \otimes \nu).
   \label{rec_dd}
\end{eqnarray}
Here we take the convention to use a highest weight to denote the
corresponding irreducible highest weight representation. It should
be clear by the notations of highest weight which Lie group or algebra
acts on.

We call the $\dinf$-modules appearing in the Fock space decomposition 
(\ref{decom_Fpairn}) {\em $n$-primitive}.
Denote by $\Od$ the category of 
all ${\dinf}_0$-diagonalizable, ${\dinf}_{+}$-locally finite 
$\dinf$-modules, with $n$-primitive
$\dinf$-modules for every $n \in \Z_{+} $
as all irreducible objects
and such that any module in $\Od$ has a Jordan-Holder composition
series in terms of $n\/$-primitive $\dinf$-modules
$( n \in \Z_{+})$. When $n = 0$ we have the trivial representation of
$\dinf$ only. 

Denote by $\Odfn$ the category of all representations
of $ O(2n)$ which decomposes into a sum of finite dimensional irreducibles.
Denote by $\Odf$ a direct sum of the categories $\Odfn$ 
for all $n \geq 0$, namely a category whose objects consist of
a direct sum of representations of $O(2n)$ in $\Odfn$, $n \geq 0$.
We introduce the following tensor product $\bigodot$ on the category $\Odf$:
given a module $U \in \Odfm$ and a module $V \in \Odfn$,
let $U \bigotimes V$ be the outer tensor product of $U$ and $V$
which becomes a $O(2m) \times O(2n)$-module. We define 
$$
  U \bigodot  V
   = \left( ind_{O(2m) \times O(2n)}^{ O(m+n)} U \bigotimes V
     \right)^{l.f.}
$$
where $ind_H^G W$ denotes the induced $G$-module from a module $W$ of
$H \subset G$ consisting of all continuous functions 
$f: G \longrightarrow W$ satisfying
$f (gh) = h^{-1} f(g),$ $ h \in H,$ $ g \in G,$
and $X^{l.f.}$ means the subset of
locally finite vectors of $X$ which are by definition 
the vectors lying in some finite dimensional $O(2m +2n)$-submodule of $X$.
One can prove the following theorem by use of Frobenius reciprocity
and the reciprocity law (\ref{rec_dd}), cf. Theorem 3.3 of \cite{W}.
\begin{theorem}
  $(\Odf, \bigodot)$ and $( \Od, \bigotimes)$ are semisimple abelian tensor
 categories. Furthermore we have an equivalence of tensor categories between
 $(\Odf, \bigodot)$ and $( \Od, \bigotimes)$ by sending
 $ V(O(2n),\nu)$ to $ L_N \left(\dinf, \Lambda^{\frak{dd}} (\nu), n \right) $.
\end{theorem}

Let us consider a second example. Following Theorem \ref{th_Fpairbd} we have
 \begin{eqnarray}
   {\cal F}^{\bigotimes (n + 1/2) }
         =  \bigoplus_{\lambda \in \Sigma(B_{n}) } 
          V(O(2n +1); \lambda) \otimes L \left(\dinf;  
                                 \Lambda^\frak{bd}_{+}(\lambda), n +1/2
                               \right).
   \label{decom_Fpairnhalf}
 \end{eqnarray}
Hence we have a see-saw pair 
\begin{eqnarray*}
 \left\{ \begin{array}{cc}
  \Bigl( O(2m +2n +1), 
    & 
    \dinf |_{c=m +n + 1/2} \Bigl) \nonumber \\
 \uparrow                          & \downarrow \nonumber\\
 \Bigl( O(2m)\times O(2n +1), 
     & \dinf |_{c=m } \bigoplus \dinf |_{c= n +1/2 } \Bigl)
                       \end{array} \right. 
\end{eqnarray*}
acting on ${\cal F}^{\bigotimes (m +n + 1/2) }$.

We call the $\dinf$-modules appearing in the Fock space decomposition 
(\ref{decom_Fpairnhalf}) {\em $(n + \hf)$-primitive}.
Denote by ${}^{\frak {bd}}{\cal O}$ the category of 
all ${\dinf}_0$-diagonalizable, ${\dinf}_{+}$-locally finite 
$\dinf$-modules, with $n$-primitive and $(n+ \hf)$-primitive
$\dinf$-modules for every $n \in \Z_{+} $
as all irreducible objects
and such that any module in ${}^{\frak {bd}}{\cal O}$ 
has a Jordan-Holder composition
series in terms of $n\/$-primitive and $(n+ \hf)$-primitive $\dinf$-modules
$( n \in \Z_{+})$. Similarly we define
a category ${}^{\frak {bd}}{\cal O}_f^k$ consisting
of all representations of $ O(k)$ which decomposes into 
a sum of finite dimensional irreducibles.
Define category ${}^{\frak {bd}}{\cal O}_f$ to be the 
direct sum of ${}^{\frak {bd}}{\cal O}_f^k$ for all $k \geq 0$.
Again we can define a tensor product on the category 
${}^{\frak {bd}}{\cal O}_f$ by the induction functor and then taking
the local finite part. Then again we have an equivalence of 
tensor categories between ${}^{\frak {bd}}{\cal O}$
and ${}^{\frak {bd}}{\cal O}_f$.

One can obtain many other reciprocity laws associated to other
dual pairs similarly.

\end{document}